\documentclass[natbib]{elsarticle}
\usepackage{setspace}
\usepackage[utf8]{inputenc}
\usepackage{color}
\usepackage{graphicx}
\usepackage{stmaryrd}
\usepackage{amsmath}
\usepackage{amssymb}
\usepackage{amsfonts}
\usepackage[ruled,vlined]{algorithm2e}
\usepackage{array}
\newcolumntype{P}[1]{>{\centering\arraybackslash}p{#1}}
\newcolumntype{M}[1]{>{\centering\arraybackslash}m{#1}}
\usepackage{tabularx}
\usepackage{stmaryrd}
\usepackage{soul}
\usepackage{diagbox}

\usepackage{hyperref}

\newcommand{\norm}[1]{\left\lVert#1\right\rVert}
\usepackage{siunitx}
\usepackage[margin=1in]{geometry}

\usepackage[utf8]{inputenc}
\usepackage{tikz}
\usepackage{tikz-dimline}
\usepackage{pgfplots}
\DeclareUnicodeCharacter{2212}{−}
\usepgfplotslibrary{groupplots,dateplot}
\usetikzlibrary{patterns,shapes.arrows}
\pgfplotsset{compat=newest}

\begin{document}
\begin{frontmatter}
	\title{Three-dimensional topology optimization of heat exchangers with the level-set method.}
	\author[LLNL]{Miguel A. Salazar de Troya}
	\author[LLNL]{Daniel A. Tortorelli}
	\author[LLNL]{Julian Andrej}
	\author[LLNL]{Victor A. Beck\corref{correspondingauthor}}
	\address[LLNL]{Lawrence Livermore National Laboratory, Livemore, CA, USA, 94550}
	\cortext[correspondingauthor]{Corresponding author, beck33@llnl.gov. }
	\begin{abstract}
		We design heat exchangers using level-set method based topology optimization.
		The heat exchange between two fluids in separate channels is maximized while constraining the pressure drop across each channel.
		The flow is modeled by an incompressible Navier-Stokes-Brinkmann equation and the heat transfer is modeled by a convection-diffusion equation with high Peclet number.
		Each fluid region is subject to its own set of Navier-Stokes-Brinkmann equations
		where the Brinkmann term models the other fluid as solid,
		thereby preventing mixing.
		A level-set defines the interface that separates the two fluids.
		The Hamilton-Jacobi equation advects the level-set, allowing for topological changes of the channels.
		The velocity of the Hamilton-Jacobi equation is defined by the shape derivatives of the cost and constraint functions with respect to normal interface perturbations.
		We present results in three-dimensions with different heat exchanger configurations and operating conditions.
	\end{abstract}
\end{frontmatter}

\newcommand{\uu}{\mathbf{u}}
\newcommand{\velocity}{\boldsymbol \xi}
\newcommand{\vv}{\mathbf{v}}
\newcommand{\xx}{\mathbf{x}}
\newcommand{\rojo}[1]{\hl{#1}}
\vspace{5mm}

\section{Introduction}

Energy, manufacturing and chemical industries use heat exchangers to efficiently transfer heat between two or more fluids traveling in separate channels.
Heat transfer occurs mostly through conduction across the channel walls, therefore maximizing the channel surface area is key for performance.
However, the pressure drop between the channel inlets and outlets increases with surface area and hence, so do pump operating costs.

Current manufacturing limitations restrict heat exchanger designs to a low dimensional design space.
The channel configuration is typically fixed and only a parametric optimization on channel cross sections or lengths can be performed
\citep{Patel2019, dixit2015review, awais2018heat, foli2006optimization}.
However, engineers now have enlarged design space thanks to the advent of additive manufacturing.
New design tools must be developed to leverage this added design freedom.

One such design tool is topology optimization.
The first applications of topology optimization were in structural mechanics wherein a given amount of material is distributed to design the stiffest structure.
Two topology optimization approaches have been used to represent the structures geometry in a mathematically tractable way.
Notably, \cite{bendsoe} formulated a material distribution problem
wherein a volume fraction field models the solid, void and intermediate material regions.
An alternative approach based on the level-set method \citep{wang2003level,allaire, SETHIAN2000489} avoids the non-physical intermediate region by using the zero isocontour of a level-set function to define the solid-void interface.
However, the ability to nucleate holes with the latter method is not straightforward, albeit it can merge existing holes.

Topology optimization has subsequently found application to design systems beyond structural mechanics, these include fluidics \citep{borrvall}, \citep{gersborg2005topology}, electromagnetics \citep{kiziltas2003topology},  \citep{zhou2010level} and acoustics \citep{wadbro}.
An extensive review of topology optimization applications can be found in \cite{deaton2014survey} and \cite{van2013level},
the latter is specific to level-set topology optimization methods.
For a more focused review on fluidic-based design, see \cite{alexandersen2020review}.

Applications of topology optimization to heat transfer systems initially considered heat conduction \citep{gersborg2006topology, li1999shape}.
However, efficient heat transfer systems rely on forced convection, caused by pressure differences \citep{dede2009multiphysics}, \citep{yoon2010topological},
and natural convection, caused by buoyancy forces generated by temperature gradients in the fluid \citep{alexandersen2014topology}.
An extensive review of the design of heat transfer systems can be found in \cite{heatreview}.

Design of heat exchangers with topology optimization has attracted much recent attention.
Although the first applications designed one flow channel \citep{haertel2017fully, kobayashi2019}, recent work designs the entire heat exchanger.
To the best of our knowledge, the first effort to design the entire system is described in the \cite{panas} thesis.
In this work, two Stokes-Brinkmann equations, one for the hot fluid channel and one for the cold water channel, are solved to model the fluids in their separate channels.
The hot and cold channel regions are complementary to each other.
As such, the channel walls separating the two fluids channels are assumed to be infinitesimally thin.
\cite{saviers2019design} designs heat exchanger in the three-dimensional setting, however, very few details of the methodology are explained due to proprietary reasons.
\cite{Kobayashi2020} extends the \cite{panas} method to three-dimensions and uses the Navier-Stokes equation to model the flow.
\cite{Hoghoj} also designs three-dimensional heat exchangers with the Navier-Stokes equation and incorporates channel wall thickness via an erosion-dilation operation.
These studies all use the volume fraction method and hence the channel wall region is modeled as a smeared interphase.
The level-set method allows for a clear definition of the channel walls.
Such an approach was first implemented by \cite{feppon2} where a single flow equation is solved and a minimum distance constraint maintains a finite channel wall thickness.
Due to this constraint, the initial design needs continuous flow paths for both the cold and hot channels between their inlets and outlets, which fixes the channel configuration.
Moreover, defining the initial design is a nontrivial task.
On the other hand, their formulation does not use an Ersatz approach to model the solid phase with a Brinkmann term, allowing them to use general purpose preconditioners to resolve the Navier-Stokes equation.

Our work solves the heat exchanger design problem using the level-set method.
Unlike \cite{feppon2}, we use an Ersatz approach and solve two fluid equations as in \cite{panas}.
Namely, the cold and hot fluid flow is modeled over the entire domain.
To simulate the cold flow, a Brinkmann term is enforced in the hot water channel which effectively negates fluid flow in this region and vice versa.
This approach allows us to start from an arbitrary initial design and optimize the channel configuration.
In an effort to increase our research dissemination, we provide the source code used in this study for replication purposes and for the basis of future research endeavors.
Notably, the codebase can be adapted for the design of generalized two-phase systems that exchange species across an hydraulically impermeable interface such as fuel cells, batteries, membrane separators, etc \citep{Nunes:2020fg,Wang:2020jk, Park:2016bt}.

In Section \ref{sec:modeling}, we present the equations to model the heat exchanger, namely the Navier-Stokes-Brinkmann and convective heat transfer equations.
Section \ref{sec:optimization} discusses the numerical strategies to solve the topology optimization problem via the level-set method.
In Section \ref{sec:examples}, our method is applied to design heat exchangers with  different configurations and operating conditions.
Finally, Section \ref{sec:conclusion} discusses the conclusions and future research directions.
This work is an extension of the work \cite{troyaTwoDimensionalTopology2021} presented at the 14th WCCM-ECCOMAS Congress,
and shares many of the details pertaining to the level-set method implementation.

\section{Governing equations}
\label{sec:modeling}
Heat exchangers simulations combine fluid mechanics and heat transfer in a highly complex manner.
For simplicity, our fluid simulations assume steady-state and incompressibility.
We thus neglect buoyancy and body forces.
The heat exchanger consists of two fluids that are confined to separate channels.
To model their respective flows, we follow the formulation used by \cite{panas, Kobayashi2020,  Hoghoj}
and solve one set of non-dimensional Navier-Stokes equations per fluid
for the non-dimensional velocity pressure pair $(\uu_i, p_i)$, where $i=C(H)$ denotes the cold (hot) fluid
\begin{equation}
	\begin{aligned}
		\uu_i \cdot\nabla\uu_i- \frac{1}{Re} \nabla^2 \uu_i + \nabla p_i & = \mathbf{0}          & \text{in} & ~ D_i                                                                         \,, \\
		\nabla \cdot \uu_i                                               & = 0                   & \text{in} & ~ D_i                                                                         \,, \\
		\uu_i                                                            & = \uu_{i_{\text{in}}} & \text{on} & ~ \Gamma_{i_{\text{in}}}                                                      \,, \\
		\frac{1}{Re}\nabla\uu_i \cdot \mathbf{n} - p\mathbf{n}           & = 0                   & \text{on} & ~ \Gamma_{i_{\text{out}}}                    \,,                                  \\
		\uu_i                                                            & =  0                  & \text{on} & ~ \Gamma_i \,,
	\end{aligned}
	\label{eq:stokes_first}
\end{equation}
where $Re=\frac{U L}{\nu}$ is the Reynolds number ($U$ is the velocity scale;
$L$ is the characteristic length of the domain and $\nu$ is the kinematic viscosity).
For the sake of brevity, we henceforth refer to the non-dimensional velocity pressure pair as simply the velocity pressure pair.
Without the loss of generality, we assume the same viscosity for the hot and cold fluids, i.e. $\nu = \nu_C = \nu_H$.

The spatial domains $D_C$ and $D_H$ are complementary subsets of the total domain $D$.
Fluid $i$ enters the channel $D_i$ through surface $\Gamma_{i_{\text{in}}}$ with an imposed velocity $\uu_{i_{\text{in}}}$ and exits through traction-free boundary $\Gamma_{i_{\text{out}}}$.
The remaining boundary $\Gamma_i=\partial D_i \setminus \Gamma_{i_{\text{in}}} \cup \Gamma_{i_{\text{out}}}$ is subjected to the no-slip boundary condition.
Our level-set strategy assumes an infinitesimally thin interface $\Gamma_i$, i.e. channel wall, separates the fluids.
It is possible to include a finite thickness channel wall, as demonstrated in the volume fraction approach of \cite{Hoghoj}.
In the context of the level-set method with an Ersatz approach applied to linear elasticity, \cite{WANG2018553} are able to represent a finite width coating using two different level-set isocontours.
Applying this technique to multiphysics nonlinear problems is deferred for future work.

A non-dimensional advection-diffusion equation models the thermal energy transport without distinguishing the phases, i.e. the non-dimensional temperature field $T$ (simply referred to as temperature henceforth)
is computed over the entire domain $D$ by solving
\begin{equation}
	\begin{aligned}
		\mathcal{L}(T) \equiv \nabla T \cdot \uu- \frac{1}{Pe} \nabla^2 T & = 0 & \text{in} & ~ D                                             \,,                                    \\
		T                                                                 & = 0 & \text{on} & ~ \Gamma_{C_{\text{in}}}   \,,                                                         \\
		T                                                                 & = 1 & \text{on} & ~ \Gamma_{H_{\text{in}}}   \,,                                                         \\
		\nabla T \cdot \mathbf{n}                                         & = 0 & \text{on} & ~ \partial D\setminus\left(\Gamma_{H_{\text{in}}}\cup\Gamma_{C_{\text{in}}} \right)\,,
	\end{aligned}
	\label{eq:temperature}
\end{equation}
where $Pe$ is the Peclet number and $\uu = \uu_C$ in $D_C$ and $\uu=\uu_H$ in $D_H$.
Again, without lose of generality, we assume the same thermal properties for the hot
and cold fluids. Adiabatic boundary conditions are enforced over all surfaces with the
exception of the channel inlets, over which the temperature is prescribed.

The cold and hot fluid domains, $D_C$ and $D_H$ respectively, are not known a priori
because their geometries evolve during the optimization.
To model their respective flow regions, we follow an Ersatz approach and simulate the flow of each fluid over the entire domain $D$ by subjecting the complementary region to a Brinkmann penalization term $\frac{1}{Da} \chi_i \uu_i$ \citep{borrvall} where $Da$ is the Darcy number, $\chi_i:D \rightarrow \{0,1\}$ are the indicator functions defined such that
\begin{align}
	\chi_C(\xx)	=
	\begin{cases}
		0 ~ \text{for}~ \mathbf{x} \in D_C \\
		1 ~ \text{for}~ \mathbf{x} \in D_H = D \setminus D_C
	\end{cases}
\end{align}
and
\begin{align}
	\chi_H(\xx) = 1 - \chi_C(\xx) \,,
	\label{eq:hotfluid}
\end{align}
As $Da\to0$, we obtain zero flow velocity in the complementary region, but excessively small $Da$ values cause numerical issues.
In practice, we assign a ``small enough'' value i.e. $0< Da \ll 1$.
Thus, we approximate the original problem in Equation \eqref{eq:stokes_first} as finding
$(\uu_i, p_i)$ such that
\newcommand{\uuin}{\uu_{i_{\text{in}}}}
\newcommand{\uuc}{\uu_{C_{\text{in}}}}
\newcommand{\uuh}{\uu_{H_{\text{in}}}}
\begin{align}
	\mathcal{L}(\uu_i, {p}_i) \equiv \uu_i \cdot{\nabla}\uu_i- \frac{1}{Re} {\nabla}^2 \uu_i + \nabla {p}_i + \frac{1}{Da}\chi_i \uu_i & = \mathbf{0} & \text{ in} & ~ D                                                                         \,, \label{eq:momentum} \\
	\nabla \cdot \uu_i                                                                                                                 & = 0          & \text{in}  & ~ D          \,, \label{eq:mass}                                                                    \\
	\uu_i                                                                                                                              & = \uuin      & \text{on}  & ~ \Gamma_{i_{\text{in}}}     \,,                         \nonumber                                  \\
	\frac{1}{Re}\nabla\uu_i \cdot \mathbf{n} - p \mathbf{n}                                                                            & = 0          & \text{on}  & ~ \Gamma_{i_{\text{out}}}                    \,,                                                    \\
	\uu_i                                                                                                                              & =  0         & \text{on}  & ~ \partial D \setminus (\Gamma_{i_{\text{in}}} \cup \Gamma_{i_{\text{out}}})  \,.\nonumber
\end{align}

We employ the Galerkin Least-Squares (GLS) stabilization method and the finite element method over the triangular partition $\mathcal{T}_h$ of $D$ to discretize and solve the Navier-Stokes equations.
As such, we find $\left(\uu_i, p_i\right) \in \mathbf{V}_i(\uuin) \times Q_i$ such that
\begin{equation}
	\begin{aligned}
		a(\chi_i; \uu_i, \vv_i)  + b(p_i, \vv_i) +  c((\uu_i, {p}_i), (\vv_i, q_i)) & =  0 \,, \\
		b(q_i, \uu_i)                                                               & = 0 \,,
	\end{aligned}
	\label{eq:stokes}
\end{equation}
for all test functions $\left(\vv_i, q_i\right) \in \mathbf{V}_{i}(0) \times Q_i$
where
\begin{equation}
	\begin{aligned}
		a(\chi_i; \uu_i, \vv_i) & = \int_{D}\left( \left(\uu_i \cdot {\nabla} \uu_i \right) \cdot \vv_i
		+ \frac{1}{Re} {\nabla} \uu_i	: {\nabla} \vv_i
		+ \frac{1}{Da} \chi_i \uu_i \cdot \vv_i\right) ~dV\,,                                           \\
		b(p_i, \vv_i)           & = \int_D p_i \nabla \cdot \vv_i ~dV \,,                               \\
	\end{aligned}
	\label{eq:stokes_expanded}
\end{equation}
and
\begin{equation}
	\begin{aligned}
		c((\uu_i, {p}_i), (\vv_i, q_i)) & = \int_D \delta_{NS_i}\mathcal{L}(\uu_i, {p}_i) \cdot \mathcal{L}(\vv_i, q_i) ~dV\,,
	\end{aligned}
\end{equation}
with
\begin{equation}
	\begin{aligned}
		\delta_{NS_i} = \beta_{GLS} \left( \frac{4 \uu_i \cdot \uu_i}{h^2} + \left(9 \frac{4 \nu}{h^2} \right)^2 + \left(\frac{1}{Da} \chi_i \right)^2\right)^{-0.5}\,.
	\end{aligned}
\end{equation}
where $h$ is the element size, and $\beta_{GLS}=0.9$ in our examples.
The finite element function spaces are
\begin{align}
	\mathbf{V}_i(\uuin) & =\left\{\mathbf{v} \in \left[H^{1}(D)\right]^{d}: \left.\mathbf{v}_{h}\right|_{\mathcal{K}} \in \left[P^{1}(\mathcal{K})\right]^{d}~\forall \mathcal{K} \in \mathcal{T}_{h} ~|~ \mathbf{v}=\mathbf{u}_{i_{\text{in}}} \text { on } \Gamma_{i_{\text{in}}}
	\text{ and }  \mathbf{v}=0 \text { on } \partial D \setminus \left(\Gamma_{i_{\text{in}}} \cup \Gamma_{i_{\text {out}}}\right)   \right\}\,,                                                                                                                                    \\
	Q_i                 & =\left\{q \in H^{1}(D): \left.q_{h}\right|_{\mathcal{K}} \in P^{1}(\mathcal{K})~\forall \mathcal{K} \in \mathcal{T}_{h} \right\}\,,
\end{align}
where $d=2,3$ is the spatial dimension.

The advection-diffusion equation, cf. Equation \eqref{eq:temperature}
is discretized in a similar manner, i.e. we find
$T\in W(T_H, T_C)=\{v_{h} \in H^{1}(D):\left.v_{h}\right|_{\mathcal{K}}
	\in P^{1}(\mathcal{K})~ \forall \mathcal{K} \in \mathcal{T}_{h} ~|~ T=T_H \text{ on }
	\Gamma_{H_{\text{in}}},~T=T_C \text{ on } \Gamma_{C_{\text{in}}} \}$ such that
\begin{align}
	a_T(\uu; T, w) + c_T(\uu; T, w)= 0\,,
\end{align}
for all $w\in W(0, 0)$, where
\newcommand{\nn}{\mathbf{n}}
\begin{equation}
	\begin{aligned}
		a_T(\uu; T, w) & = \int_{D} \left(\uu \cdot {\nabla} T w + \frac{1}{Pe} {\nabla} T
		: {\nabla} w \right) ~dV \,,                                                       \\
		\label{eq:temp_gls}
	\end{aligned}
\end{equation}
and
\begin{equation}
	\begin{aligned}
		c_T(\uu; T, w) & = \int_D \delta_{T}\mathcal{L}_T(T) \cdot \mathcal{L}_T(w) ~dV\,,
	\end{aligned}
\end{equation}
with
\begin{equation}
	\begin{aligned}
		\delta_{T} = \beta_{GLS} \left( \frac{4 \uu \cdot \uu}{h^2} + \left(9 \frac{4 }{h^2Pe} \right)^2 \right)^{-0.5}
	\end{aligned}
\end{equation}
wherein $\uu = \uu_H + \uu_C$.

\section{Topology optimization of heat exchangers}
\label{sec:optimization}
Our heat exchanger is designed to maximize heat transfer from hot to the cold fluid.
To do this, we solve a topology optimization problem to determine the cold fluid domain $D_C$, and hence the complementary hot fluid domain $D_H$.

\subsection{Shape optimization problem}
We formulate the topology optimization problem as a shape optimization problem.
Letting $\mathcal{O}$ denote the set of all admissible domains, we find a $D_C\in\mathcal{O}$ that solves
\newcommand{\VV}[1]{\mathbf{V}_{#1}}
\begin{subequations}
	\label{eq:opti_problem_ls}
	\begin{flalign}
		\underset{D_C \in \mathcal{O}}{\text{max}}~ &J(D_C)  = \int_{\Gamma_{C_{\text{out}}}} T \uu_C \cdot \mathbf{n}~dA\,,&&  \label{eq:cost_func} \\
		\text{s.t.}& ~(\uu_C, \uu_H, p_C, p_H, T)  \in \VV{C}(\uuc)\times \VV{H}(\uuh) \times Q_C \times Q_H \times W(T_H, T_C) ~\text{satisfy} ~  &&  \nonumber\\
		&a(\chi_C; \uu_C, \vv_C) + b(\vv_C, p_C) + b(\uu_C, q_C) + c((\uu_C, {p}_C), (\vv_C, q_C))  = 0  \label{eq:stokes1_ls}  \,,&&   \\
		&a(1 - \chi_C; \uu_H, \vv_H) + b(\vv_H, p_H) + b(\uu_H, q_H)+ c((\uu_H, {p}_H), (\vv_H, q_H))  = 0  \label{eq:stokes2_ls} \,,&& \\
		&a_T(\uu_C + \uu_H; T, w) + c_T(\uu_C + \uu_H; T, w)   = 0 \label{eq:temp_ls}\,, &&  \\
		\text{for all}& ~(\vv_C, \vv_H, q_C, q_H, w ) \in \VV{C}(\mathbf{0}) \times \VV{H} (\mathbf{0})\times Q_C \times Q_H \times W(0, 0) \nonumber ~\text{and}&&  \\
		&G_1(D_C)  = \int_{\Gamma_{C_{\text{in}}}} p_C ~dA  - \int_{\Gamma_{C_{\text{out}}}} p_C ~dA   \leq P_{\text{drop}}\,,&&\label{eq:constraint_1}\\
		&G_2(D_C) = \int_{\Gamma_{H_{\text{in}}}} p_H ~dA  - \int_{\Gamma_{H_{\text{out}}}} p_H ~dA  \leq P_{\text{drop}}\label{eq:constraint_2}\,,&&
	\end{flalign}
\end{subequations}
As seen above, our goal is to maximize the heat flux accross the cold fluid outlet
$\Gamma_{C_{\text{out}}}$ subject to the constraints $G_1(D_C) \leq P_{\text{drop}}$ and $G_2(D_C) \leq P_{\text{drop}}$
that limit the pressure drop across each fluid channel to be less than $P_{\text{drop}}$.
These constraints prevent designs with excessive channel surface area and thereby regularize the optimization problem.
The shape optimization is performed with respect to the cold fluid domain $D_C$ since that necessarily defines the hot fluid domain $D_H=D\setminus D_C$.

We use first-order methods to solve Equation \eqref{eq:opti_problem_ls}.
The analytical shape sensitivities $DJ(D_C), DG_1(D_C)$ and $DG_2(D_C)$ are obtained by differentiating the continuum equations with respect to the design and then discretizing them for evaluation purposes (differentiate-discretize).
An alternative approach, commonly used in the topology optimization community, discretizes the PDEs first and then differentiates the discretized equations with respect to the design parameters \citep{kreissl2011explicit} (discretize-differentiate).
Our choice allows us to use the automatic shape derivative calculation
in the Unified Form Language (UFL) \citep{ham2019automated} within
the pyadjoint library \citep{dokken2020automatic}, \citep{Mitusch2019}.

We use a reduced space sensitivity formulation and hence we consider the dependency of the primal response $\uu_C, p_C, \uu_H, p_H$ and $T$ on the design $D_C$ so that, e.g.
$J(D_C) = \hat{J}(D_C, \uu_C(D_C), p_C(D_C), \uu_H(D_C), p_H(D_C), T(D_C))$.
Following \cite{laurain2018level}, a general functional $J(D_C)$ evaluated on the
design $D_C$ admits a shape derivative $DJ(D_C)[\velocity]$ in the direction
$\velocity \in \Psi(D)= \{ \velocity \in C^s(\mathbb{R}^d)~ |~ \velocity \cdot
	\mathbf{n}_{|\partial D} = 0 \}$ (with integer $s \geq 0$) in the form
\begin{equation}
	D J(D_C)[ \velocity]=\int_{D} \left(\Pi: \nabla \velocity+\iota \cdot \velocity\right) ~dV \,,
	\label{eq:shapeder}
\end{equation}
where $\Pi$ and $\iota$ are sufficiently smooth tensor and vector fields.
Given enough regularity assumptions, the Hadamard-Zol\'esio theorem
\citep{delfour2011shapes} shows that Equation \eqref{eq:shapeder} only depends
on the restriction to the boundary
\newcommand{\zetav}{\boldsymbol{\zeta}}
\begin{equation}
	DJ(D_C)[\velocity]=\int_{\Gamma} \llbracket \Pi \rrbracket \velocity \cdot n ~dA \,.
	\label{eq:shapederbound}
\end{equation}
We opt however, to use the volumetric Equation \eqref{eq:shapeder} because it is easier to compute with pyadjoint.
Moreover, it also offers better accuracy \citep{choi1986design, hiptmair2015comparison}, as it is consistent with the discretize and then differentiate approach \citep{berggren2010unified}.


The domain shape derivative $DJ(D_C)[\velocity]$ in Equation \eqref{eq:shapeder}
requires the specification of $\zetav_J$, i.e. the shape velocity field throughout $D$,
rather than just over the interface $\Gamma$ as required by Equation \eqref{eq:shapederbound}.
This domain specification of $\zetav_J$ is known as the extension problem \citep{van2013level}.
Here, we follow \cite{gournay2006velocity} and use a Hilbertian extension which helps accelerate convergence of the optimization problem, regularizes the shape derivative and endows it with a scalar product.
Once we formulate the shape derivative $DJ(D_C)[\velocity]$, we proceed to
find the velocity field $\zetav_J \in \Psi_h(D)$ such that
\begin{equation}
	\zetav_J = \underset{\velocity \in \Psi_h(D)}{\text{argmin}} DJ(D_C)[\velocity]\,,
	\label{eq:optimal_velocity}
\end{equation}
where $\Psi_h(D)= \{ \velocity \in \left[H^1(D)\right]^d\}$.
The solution to \eqref{eq:optimal_velocity} is obtained by finding $\zetav_J\in \Psi_h(D)$ such that
\begin{equation}
	\begin{aligned}
		b_{\zetav}(\zetav_J, \velocity) & = DJ(D_C)[\velocity]               ~~ \forall ~\velocity \in \Psi_h(D) \,,
	\end{aligned}
	\label{eq:regularization}
\end{equation}
where
\begin{equation}
	b_{\zetav}(\zetav, \velocity) = \int_{D} \gamma \nabla \zetav: \nabla \velocity ~dV + \int_{D} \boldsymbol\zetav \cdot \velocity ~dV + c_1 \int_{\partial D} \left(\zetav \cdot \mathbf{n} \right) \left(\velocity \cdot \mathbf{n} \right) ~dA\,.
	\label{eq:inner}
\end{equation}
The parameter $\gamma$ regularizes the problem and helps to accelerate convergence of the optimization \citep{gournay2006velocity,burger2003framework}.
It must be large enough so that the support of $\zetav$ extends beyond the interface $\partial D_C$, but small enough so that values of $\zetav$ associated with neighboring interface regions do not interfere with each other.
Note that the Dirichlet boundary condition $\velocity\cdot\mathbf{n}_{|\partial D} = 0$ is imposed by means of a penalty method where $c_1 \gg 1$.
The same regularization is performed on $DG_1(D_C)$ and $DG_2(D_C)$
to obtain $\zetav_{G_1}$ and $\zetav_{G_2}$.

Once $\zetav_{J}, \zetav_{G_1}$ and $\zetav_{G_2}$ are calculated, it is necessary to calculate a search direction $\boldsymbol{\theta} \in \Psi_h(D)$ that decreases the cost function and respects the constraint inequalities $G_1(D_C)\leq P_{\text{drop}}$ and $G_2(D_C)\leq P_{\text{drop}}$ in Problem \eqref{eq:opti_problem_ls}.
It is common in the literature to minimize a Lagrangian obtained by augmenting the cost function with the constraint functions weighted by Lagrange multiplier estimates, e.g. $\lambda_1^{\text{est}}\left(G_1(D_C) - P_{\text{drop}} \right)$.
This approach does not ensure satisfaction of the constraints, but it is sometimes sufficient to obtain reasonable designs.
Other approaches use the Augmented Lagrange method \citep{allaire2014multi}, the sequential linear programming method \citep{dunning2015introducing} and the null space gradient flow algorithm \citep{feppon:hal-01972915}.
In this work, we opt for the last option, cf. \cite{feppon:hal-01972915} for details.

\subsection{Level-set method}
In our level-set approach to topology optimization, we use the differentiable level-set function $\phi:D \rightarrow \mathbb{R}$ to define $D_C$ and $\Gamma$ such that
\begin{align}
	D_C    & = \left\{\mathbf{x} \in D ~|~ \phi(\xx) < 0 \right\} \nonumber \,, \\
	\Gamma & = \left\{\mathbf{x} \in D ~|~ \phi(\xx) = 0 \right\}\nonumber  \,.
\end{align}
Since $D_H$ is complementary of $D_C$, we have $D_H= \left\{\mathbf{x} \in D ~|~ \phi(\xx) > 0 \right\}$.
As seen above, the channel wall, i.e interface $\Gamma$ is defined by the $\phi(\xx)=0$ level-set.
Changes in $\phi$ thusly change $D_C$ and $\Gamma$ and hence the heat exchanger.

The level-set function space $\Phi(D)$ is discretized over a regular mesh via linear interpolation, i.e. Lagrange shape functions.
Because of this, the material interface $\Gamma$ is rarely aligned with the element boundaries.
This well known dilemma is resolved by using conformal meshes that track the interface
\citep{allaire2011topology,yamasaki2011level},
applying immersed boundary techniques \citep{burman2018shape,villanueva2017cutfem,
	belytschko2003topology,kreissl2012levelset} and
by implementing an Ersatz material approach \citep{allaire2004structural}.
For more details on these three approaches, we refer to \cite{van2013level}.

In this work, we implement the Ersatz method.
We proceed by defining the indicator function $\chi_H = H\circ\phi: D \rightarrow \mathbb{R}$, where $H$ is the unit step function.
As such $\int_{D_H}  f(x) \,  dV = \int_D f(x)  \chi_H \, dV$ and $\int_{D_C}  f(x) \,  dV = \int_D f(x) \left(1 - \chi_H\right) \, dV$.
In the finite element discretization, $\phi$ is interpolated with  linear Lagrange elements.
Therefore, in elements where the level-set nodal values $\phi$ change sign, the effective penalization coefficient $\frac{1}{Da} \chi_H$ varies linearly between $\frac{1}{Da}$ and zero.
This results in a blurred interphase region rather than a sharp interface between the hot and cold fluid channels.
Fortunately, the size of the interphase region is reduced through mesh refinement.
The lack of differentiability of $\chi_H$ is only a concern in the
unlikely scenario when $\phi(\xx)=0$ at an element node.
If, e.g. all the element nodal level-set values are negative, except one which equals zero, then the design will be non-differentiable.
Indeed, a ``positive'' design change to the zero level-set valued node will move the entire element into the cold fluid domain whereas a ``negative'' change will keep the entire element in the hot domain.
The effect of this anomaly is also lessened with mesh refinement.

The level-set function $\phi_n: D \rightarrow \mathbb{R}$ at the optimization  iteration $n$ is advected in the search direction $\boldsymbol{\theta}$ for pseudo time $\hat{t}_f$ via the Hamilton-Jacobi equation,
\begin{equation}
	\begin{aligned}
		\frac{\partial \phi}{\partial \hat t} & =\hat{\boldsymbol\theta} \cdot \nabla \phi+ \frac{1}{Pe_{HJ}}\nabla^2 \phi & 0 \leq \hat{t} \leq \hat{t}_f \,, \\
		\phi(\xx, 0)                          & = \phi_n(\xx) \,,
		\label{eq:hamilton_jacobi}
	\end{aligned}
\end{equation}
such that $\phi_{n+1}(\xx)=\phi(x, \hat{t}_f)$.
In the above, we introduced the scales: $\hat{\boldsymbol\theta}=\frac{\boldsymbol\theta}{\theta_{\text{max}}}$, $\hat t = \frac{t}{L / \theta_{\text{max}}}$ and $\hat{\mathbf{x}} = \frac{\mathbf{x}}{L}$ where $\theta_{\text{max}}=\underset{D}{\text{max}} \left(\left|\theta_x\right|+\left|\theta_y\right| + \left|\theta_z\right|\right) $.
Note that the level-set $\phi$ does not require scaling.
The artificial conductivity $k$ in the Peclet number $Pe_{HJ}=\frac{ L \theta_{\text{max}}}{k}$ is included for stabilization purposes.
We solve Equation \eqref{eq:hamilton_jacobi} using the same GLS stabilization
scheme that is used to solve Equation \eqref{eq:temp_gls}. We represent Equation
\eqref{eq:hamilton_jacobi} in the operator form $\phi_{n+1} = M(\phi_n, \boldsymbol\theta, \hat{t}_f)$.
A line search is performed to determine the pseudo time $\hat{t}_f$ using the merit function $merit(\phi)$ described in \cite{feppon:hal-01972915}. $\hat{t}_f$ is accepted if $merit(\phi_{n+1}) \leq merit(\phi_{n})$, otherwise $\hat{t}_f \leftarrow \frac{\hat{t}_f}{2} $ and $\phi_{n+1} =M(\phi_n, \boldsymbol\theta, \frac{\hat{t}_f}{2})$ is reevaluated.
Line search iterations continue until $merit(\phi_{n+1}) \leq merit(\phi_{n})$ or the number of line search iterations reaches $maxtrials$.
The time integration is performed using the backward Euler scheme with adaptive time stepping that is implemented in the PETSc TS library \citep{abhyankar2018petscts}.
In future work, we will explore using more efficient time integrators from the PETSc TS library.

We apply a reinitialization of the level-set whenever $\tau \ge d_{\text{max}}$ to prevent $|\nabla \phi|$ from becoming too large or too small, where $\tau$ is a metric of the distance traveled between reinitializations, i.e. $\tau \leftarrow \tau + \theta_{\text{max}}\, \hat{t}_f$ at the end of each optimization iteration.
The user defined parameter $d_\text{max}$ controls the frequency of reinitialization.
In this work, we assign $d_\text{max}=0.08$.
The reinitialized level-set $\varphi$ solves the minimization problem
\begin{equation}
	\begin{aligned}
		\varphi = \underset{\varphi \in \Phi_{\Gamma}(D)}{\text{argmin}} \int_{D} \eta(|\nabla \phi|)~dV
		\label{eq:reinit_min}
	\end{aligned}
\end{equation}
where the potential $\eta$
\begin{equation}
	\begin{aligned}
		\eta(\alpha)=\left\{\begin{array}{ll}\frac{1}{2}(\alpha-1)^{2} & \text { if }\alpha > 1 \\ \frac{\alpha^{3}}{3}-\frac{\alpha^{2}}{2}+\frac{1}{6} & \text { if }\alpha \leq 1\end{array}\right.
	\end{aligned}
\end{equation}
includes the residual of the Eikonal equation and a special extension to avoid the singularity at $|\nabla \phi| = 0$ \citep{adamsHighorderEllipticPDE2019}.
The function space $\Phi_{\Gamma}(D)$ contains the interface condition $\varphi\mid_{\Gamma}=0$ such that it leaves $\Gamma$ unchanged.
Taking the variational of Equation \eqref{eq:reinit_min} yields the elliptic problem of finding $\varphi\in\Phi_{\Gamma}(D)$ such that
\begin{equation}
	\begin{aligned}
		\int_{D} \iota(|\nabla \varphi|) \nabla \varphi \cdot \nabla v ~dV & = 0 \\
		\label{eq:signed_distance}
	\end{aligned}
\end{equation}
for all $v \in \Phi_{\Gamma}(D)$.
The diffusion coefficient in the above Equation \eqref{eq:signed_distance} is
\begin{align}
	\iota(\alpha)=\left\{\begin{array}{ll}1-\frac{1}{ \alpha } & \text { if }\alpha >1 \\ 1-(2-\alpha) & \text { if }\alpha \leq 1\end{array}\right.
	\label{eq:reinit_elliptic}
\end{align}
The numerical procedure used to solve Equation \eqref{eq:signed_distance} is described in Appendix \ref{sec:sds}.
Notably, it is not straightforward to enforce the $\varphi\mid_{\Gamma}=0$ condition.

Finally, we summarize the optimization procedure in Algorithm \ref{alg:optimization}.

\begin{algorithm}[h!]
	\SetAlgoLined
	Initialize level-set function\;
	$n=0$\;
	\While{$ \norm{\boldsymbol{\theta}}_{\Psi_h} < tol$ or $n < maxiter$}{
		Solve Equations \eqref{eq:stokes1_ls}, \eqref{eq:stokes2_ls} and \eqref{eq:temp_ls}\;
		Evaluate \eqref{eq:cost_func}, \eqref{eq:constraint_1} and \eqref{eq:constraint_2}, i.e. $J(D_C)$, $G_1(D_C)$ and $G_2(D_C)$.\;
		Calculate the shape derivatives $DJ(D_C), DG_1(D_C)$ and $DG_2(D_C)$\;
		Regularize the shape derivatives with Equation \eqref{eq:regularization} to obtain the velocities $\zetav_J$, $\zetav_{G_1}$, $\zetav_{G_2}$\;
		Apply Null-space search algorithm to obtain the search direction $\boldsymbol{\theta}$\;
		Update level-set function $\phi_{n+1} \gets M(\phi_n, \boldsymbol{\theta}, \hat{t}_f)$ using Equation \eqref{eq:hamilton_jacobi} and the line search\;
		\For{k=1 ... \text{maxtrials}}{
			\If{\text{merit}($\phi_{n+1}$) $\leq$ \text{merit}($\phi_n$)}{
				break\;
			}
			Update level-set function $\phi_{n+1} \gets M(\phi_n, \boldsymbol{\theta},\frac{\hat{t}_f}{2^k})$\;
		}
		$\tau \leftarrow \tau + \theta_{\text{max}}\, \hat{t}_f$\;
		\If{$\tau \ge d_{\text{max}}$}{
			Reinitialize the level-set function $\phi_n$ via Equation \eqref{eq:signed_distance}\;
			$\phi_n = \varphi$\;
			$\tau = 0$\;
		}
		$n = n + 1$\;
	}
	\caption{Level-set method for topology optimization}
	\label{alg:optimization}
\end{algorithm}

\section{Numerical Examples}
\label{sec:examples}
\begin{figure}[h!]
	\centering
	\includegraphics[scale=0.2]{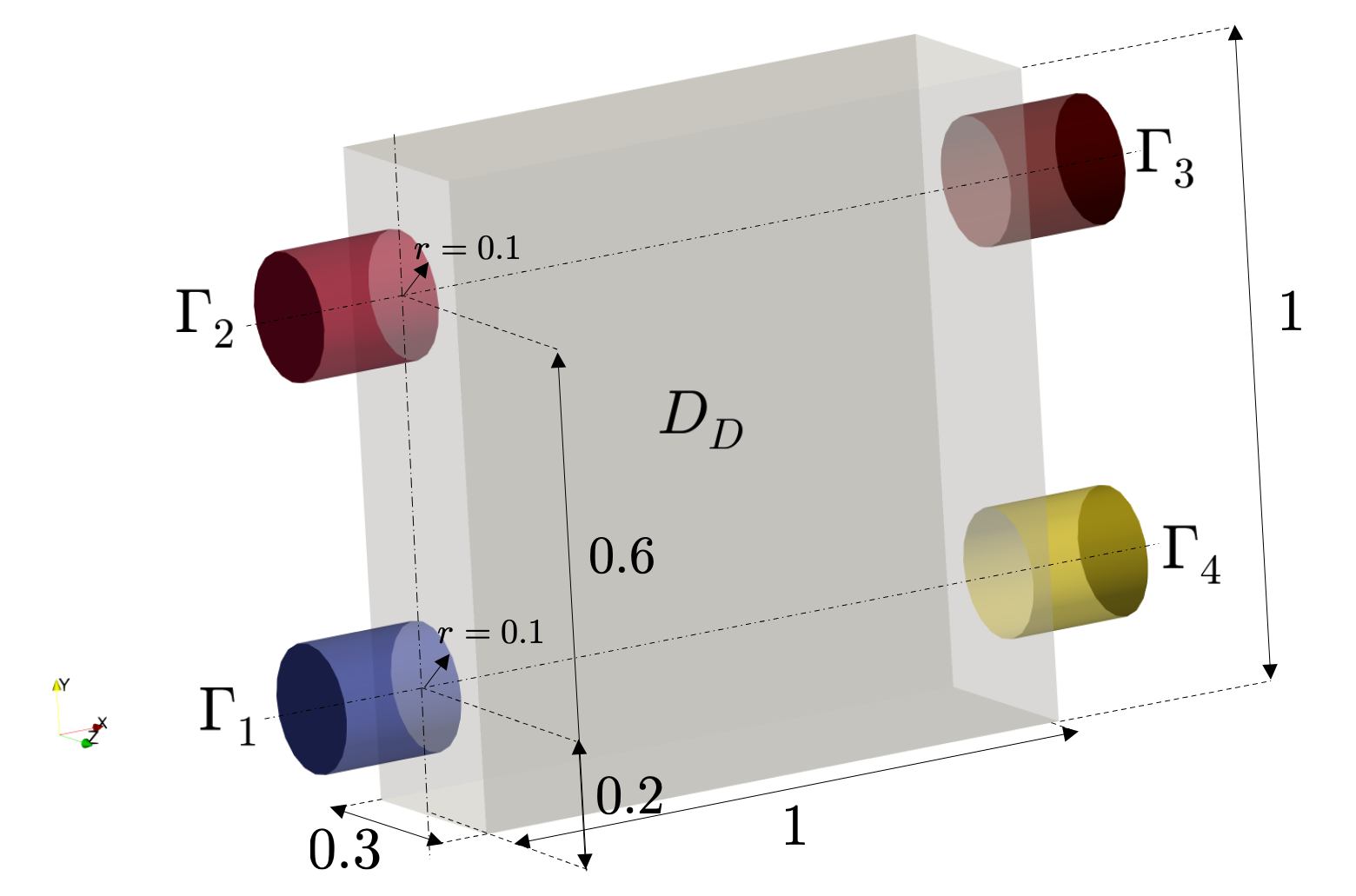}
	\caption{Design domain $D_D$ and fixed colored regions.}
	\label{fig:box_domain}
\end{figure}

\begin{figure}[h!]
	\centering
	\includegraphics[scale=0.4]{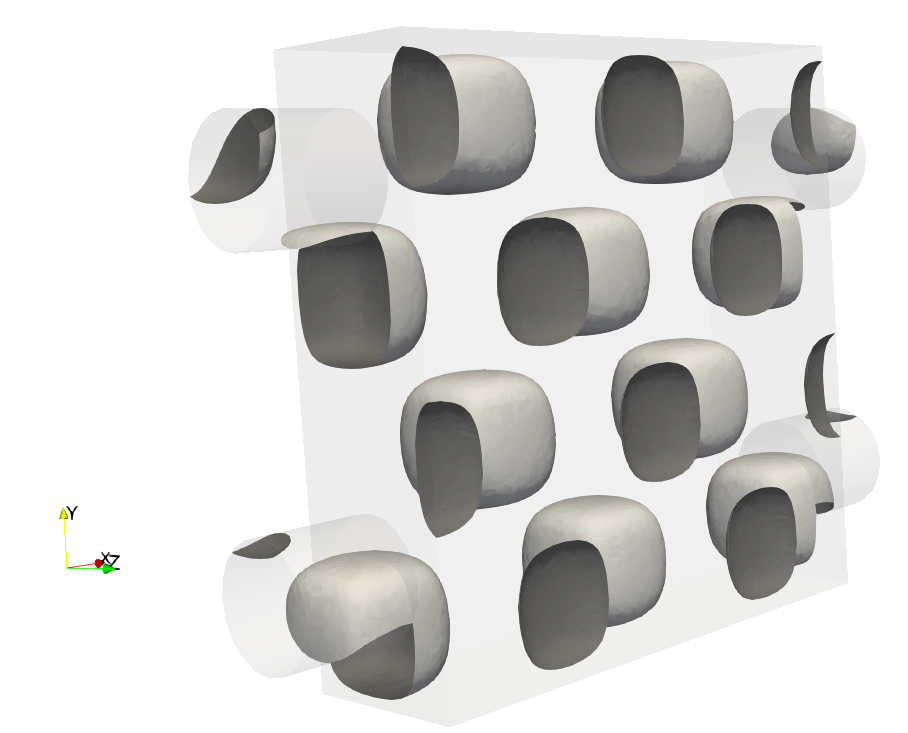}
	\caption{Initial design front view shows cold in black and hot in white.}
	\label{fig:initial_design}
\end{figure}

Our simulation domain is similar to the one found in \cite{kobayashi2019}, cf. Figure
\ref{fig:box_domain}, although we simulate the entire domain instead of assuming reflective symmetry about the XY plane.
The boundary surfaces $\Gamma_1, \Gamma_2, \Gamma_3$ and $\Gamma_4$ are assigned as cold and hot fluid inlets and outlets
($\Gamma_{C_{\text{in}}}, \Gamma_{H_{\text{in}}}, \Gamma_{C_{\text{out}}} $,
and $\Gamma_{H_{\text{out}}}$) depending on the heat exchanger configuration requested.
For the cold fluid equation, the non homogeneous Dirichlet condition $\uu_{C_{\text{in}}}$ over the inlet $\Gamma_{C_{\text{in}}}$ is a horizontal parabolic profile with a maximum non-dimensional velocity $V_{\text{max}}=1$.
At the outlet $\Gamma_{C_{\text{out}}}$, a traction free condition is imposed.
The remaining boundary $\partial D \setminus(\Gamma_{C_{\text{in}}} \cup \Gamma_{C_{\text{out}}})$ is subject to the no-slip condition $\uu_C = 0$.
The hot fluid equation has the same boundary conditions, but on $\Gamma_{H_{\text{in}}}$, $\Gamma_{H_{\text{out}}}$ and $\partial D \setminus (\Gamma_{H_{\text{in}}} \cup \Gamma_{H_{\text{out}}})$ respectively.
The Peclet number is $Pe=5\times10^3$.
The fixed colored regions $D \setminus D_D$ are used to ensure fully developed flow boundary conditions in the inlet and outlet regions.
Zero velocity is enforced in the complementary colored regions to model them as solid.
We constrain the velocity in these regions rather than applying the Brinkmann term to facilitate convergence of the linear solver.
Our initial design is $\phi(x, y, z) = \text{sin}(4 \pi y) \text{cos}(4 \pi x) \text{sin}(4 \pi z) - \frac{1}{5}$, cf. Figure \ref{fig:initial_design}, the regularization parameter is $\gamma=0.4$, and the Darcy number is $Da = 1 \times 10^{-5}$ 1/s.
We found this $Da$ value returns the most reasonable results.
Even though larger $Da$ values approximate the solid domain better, we observed that the preconditioner struggles to converge if disconnected phases are present, cf. Figure \ref{fig:initial_design}.
To alleviate this issue, we increase $Da$ by a factor of 10 after both pressure drop constraints are satisfied for the first time.
In our experience, this generally means that both channels are continuous from their inlets to their outlets and the flow is smooth.
The level-set field $\phi: D\rightarrow \mathbb{R}$ and the velocity fields $\zetav_J$, $\zetav_{G_1}$ and $\zetav_{G_2}$ components are discretized on the same mesh as the Navier-Stokes-Brinkmann equations \eqref{eq:stokes}.
The mesh is created with Gmsh \citep{geuzaine2009gmsh} and contains a total of 1,406,688 tetrahedral elements.
Our designs had mostly converged to a local optimum in less than $maxiters=1200$ iterations, cf. Algorithm \ref{alg:optimization}.

The Navier-Stokes-Brinkmann and heat transfer equations
cf. Equation (\ref{eq:stokes1_ls}, \ref{eq:stokes2_ls} and \ref{eq:temp_ls}) are discretized in Firedrake \citep{Rathgeber2016,Luporini2016,Homolya2017},
which uses PETSc \citep{petsc-efficient, petsc-user-ref,Dalcin2011} as the backend for the linear algebra.
We use the SIMPLE preconditioner \citep{patankarCalculationProcedureHeat1972} with FGMRES \citep{saadfgmres} to resolve the Navier-Stokes-Brinkmann equation
and the Hypre BoomerAMG preconditioner for the convection-diffusion equation.
The simulations were run on ten 2.60 GHz Intel XeonE5-2670 processors.

We study three different configurations for the heat exchanger: parallel, counter and U-flow whose boundary condition assignments are summarized in Table \ref{tab:configurations}.
The pressure drop constraint uses $P_{\text{drop}}=2.0$ and the operating condition uses $Re=10$.
Figures \ref{fig:parallel_clear} and \ref{fig:parallel_diffuse} show the optimized design for a parallel heat exchanger where the color indicates the fluid temperature.
The optimized design does not exhibit any domain symmetry, cf. Figure \ref{fig:parallel_side}.
Indeed, the optimizer finds a design that allows both fluids to cross each other near the XY plane.
Our next design considers a counter flow configuration, cf. Figures \ref{fig:counter_clear}-\ref{fig:counter_side}.
The design shares features from the parallel design and again, is not symmetric.
Finally, the U-flow configuration, cf. Figures \ref{fig:u_flow_clear}-\ref{fig:u_flow_side}, optimized design contains a network of pipes that cross each other near the XY symmetry plane as well.
Table \ref{tab:cost_function_values} denotes the cost function values for the three designs, where it is seen that the parallel heat exchanger design is slightly more efficient than the others.
This is unlike commercially available heat exchangers, for which counter flow designs are the most efficient.
We conjecture that this difference is due to: the design domain being more compact instead of elongated, the low Reynolds number and the dependency on the initial design.

\begin{table}[]
	\centering
	\begin{tabular}{l|llll}
		\diagbox{Configuration}{Flow surface} & Cold inlet $\Gamma_{C_{\text{in}}}$  & Hot inlet $\Gamma_{H_{\text{in}}}$ &
		Cold outlet $\Gamma_{C_{\text{out}}}$ & Hot outlet $\Gamma_{H_{\text{out}}}$                                                                \\
		\hline
		Parallel                              & $\Gamma_1$                           & $\Gamma_2$                         & $\Gamma_4$ & $\Gamma_3$ \\
		Counter                               & $\Gamma_4$                           & $\Gamma_2$                         & $\Gamma_1$ & $\Gamma_3$ \\
		U-flow                                & $\Gamma_3$                           & $\Gamma_2$                         & $\Gamma_4$ & $\Gamma_1$ \\
	\end{tabular}
	\caption{Flow boundary conditions for the heat exchanger configurations.}
	\label{tab:configurations}
\end{table}

\newcommand{\imagesize}{0.15}
\begin{figure}[h!]
	\centering
	\includegraphics[scale=\imagesize]{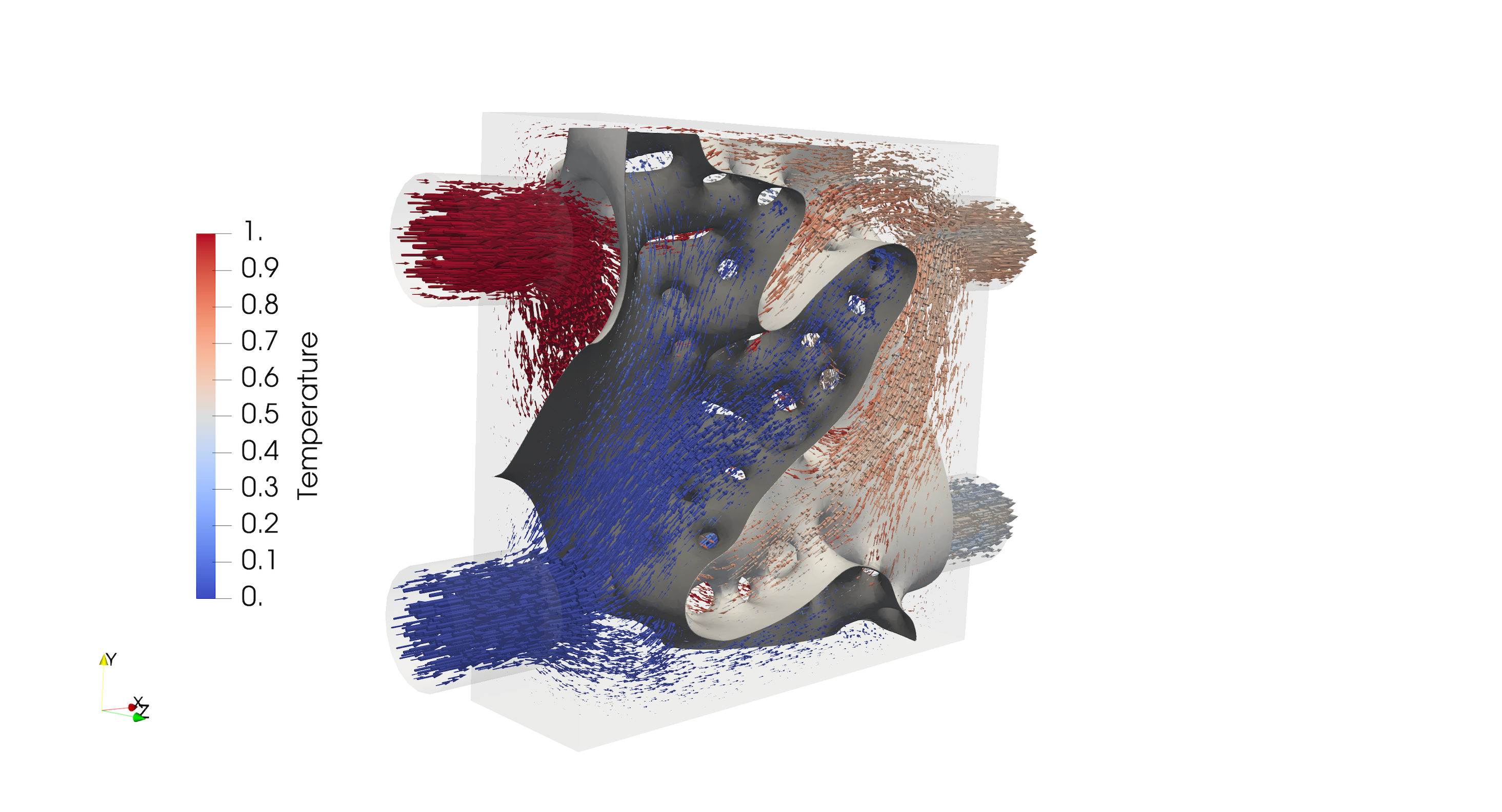}
	\caption{Optimized design, flow field and temperature for the parallel flow heat exchanger with $P_{\text{drop}} = 2.0$ and $Re=10$, front view.
		Cold membrane side in black, hot side in white.}
	\label{fig:parallel_clear}
\end{figure}
\begin{figure}[h!]
	\centering
	\includegraphics[scale=\imagesize]{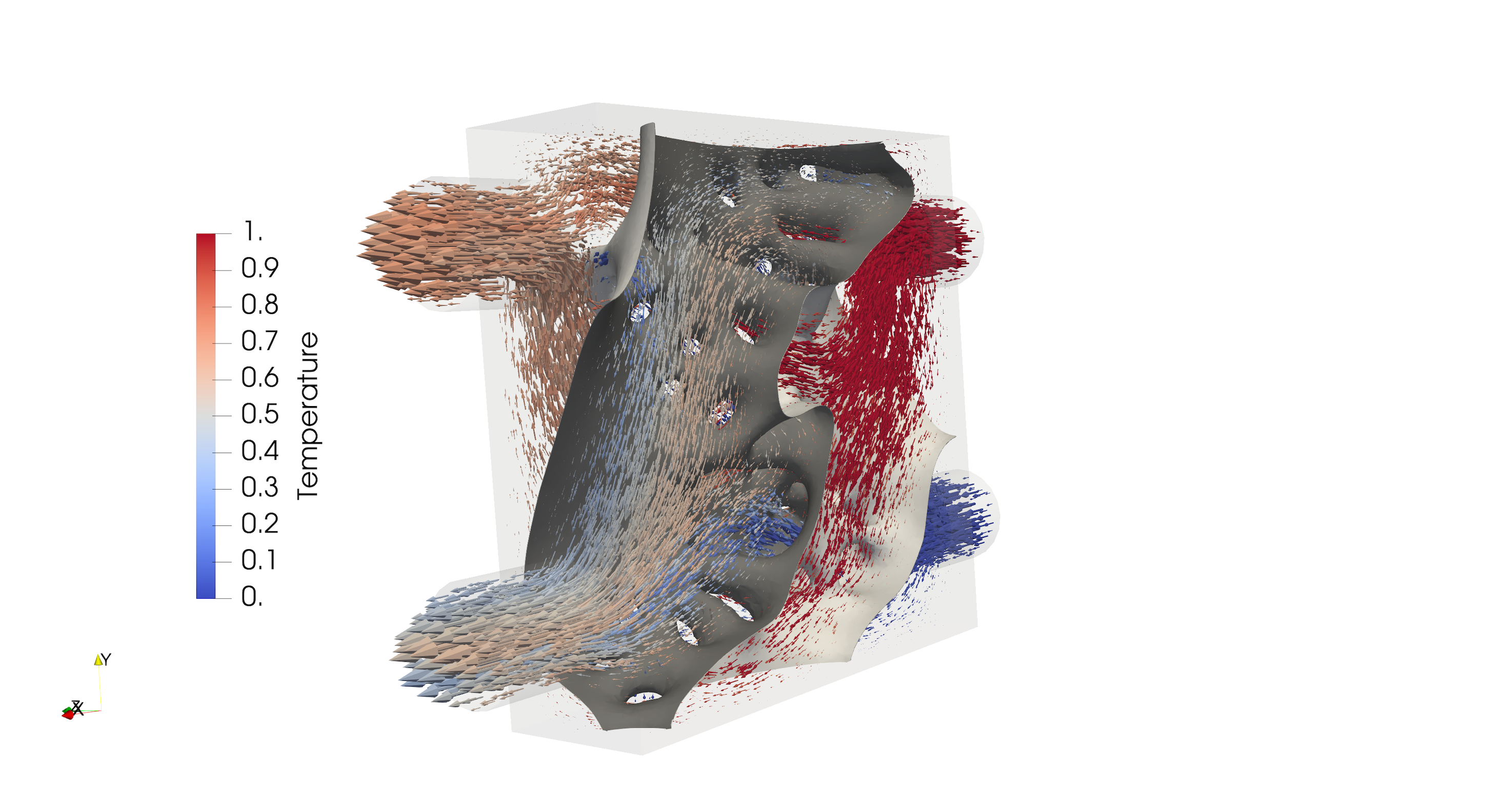}
	\caption{Optimized design, flow field and temperature for the parallel flow heat exchanger with $P_{\text{drop}} = 2.0$ and $Re=10$, back view.
		Cold membrane side in black, hot side in white.}
	\label{fig:parallel_diffuse}
\end{figure}
\begin{figure}[h!]
	\centering
	\includegraphics[scale=\imagesize]{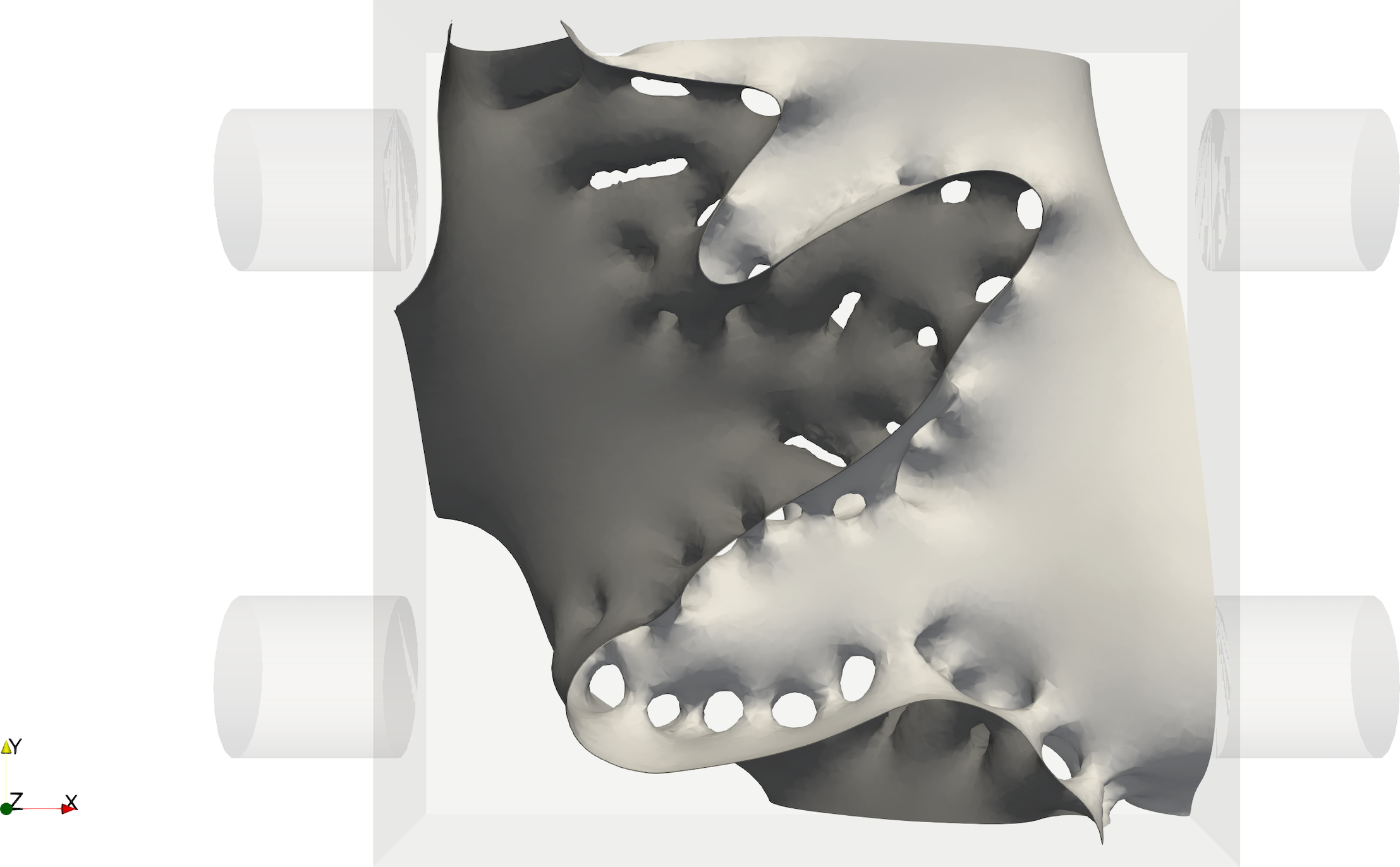}
	\caption{Optimized design for the parallel flow heat exchanger with $P_{\text{drop}} = 2.0$ and $Re=10$.
		Cold membrane side in black, hot side in white.}
	\label{fig:parallel_side}
\end{figure}

\begin{table}[]
	\centering
	\begin{tabular}{l|l|l|l}
		Configuration & Parallel & Counter & U-flow  \\
		\hline
		Cost function & -0.489   & -0.4265 & -0.4014 \\
	\end{tabular}
	\caption{Cost function values for the optimized design configurations.}
	\label{tab:cost_function_values}
\end{table}

\begin{figure}[h!]
	\centering
	\includegraphics[scale=\imagesize]{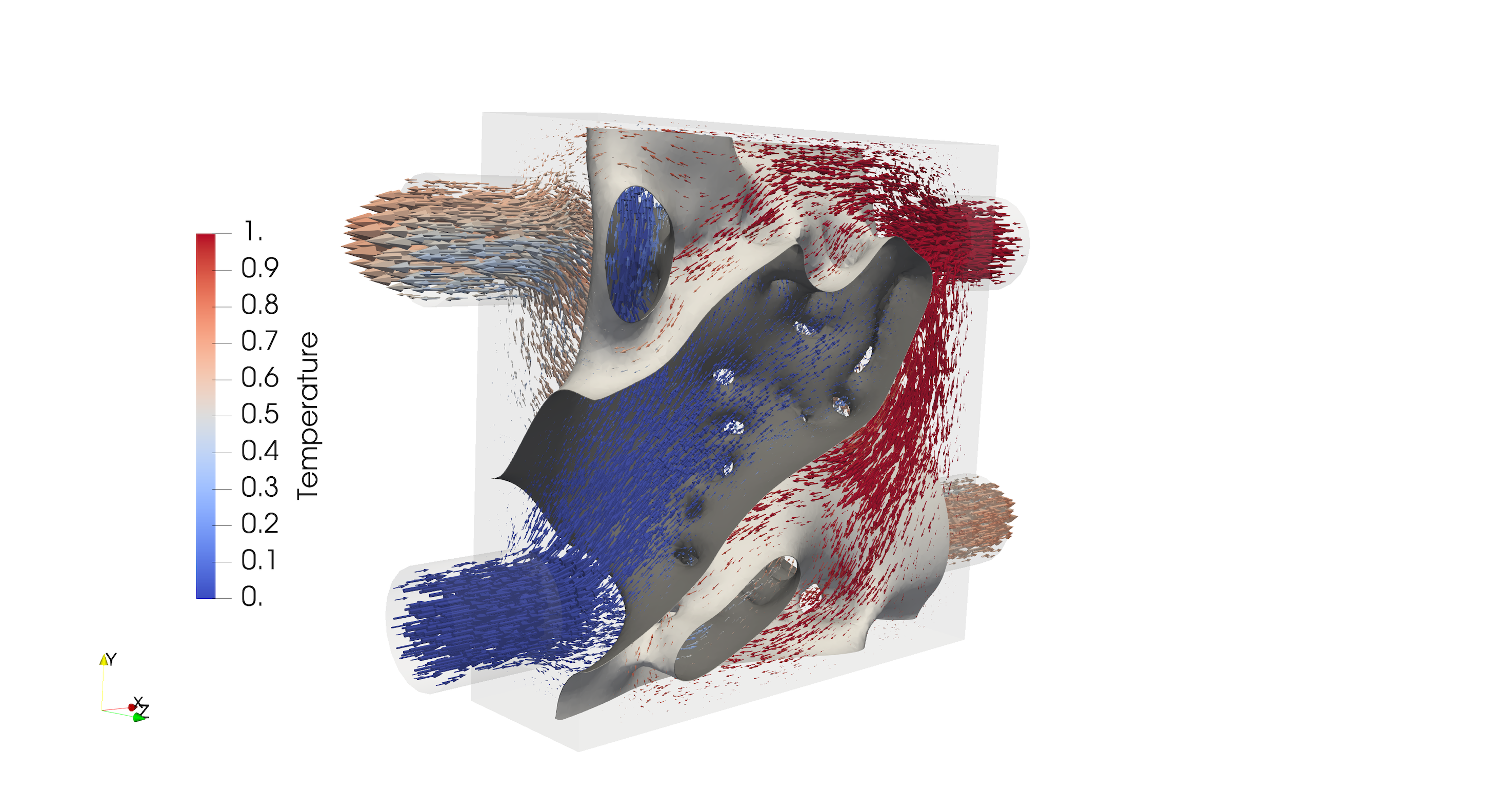}
	\caption{Optimized design, flow field and temperature for the counter flow heat exchanger with $P_{\text{drop}} = 2.0$ and $Re=10$, front view.
		Cold membrane side in black, hot side in white.}
	\label{fig:counter_clear}
\end{figure}
\begin{figure}[h!]
	\centering
	\includegraphics[scale=\imagesize]{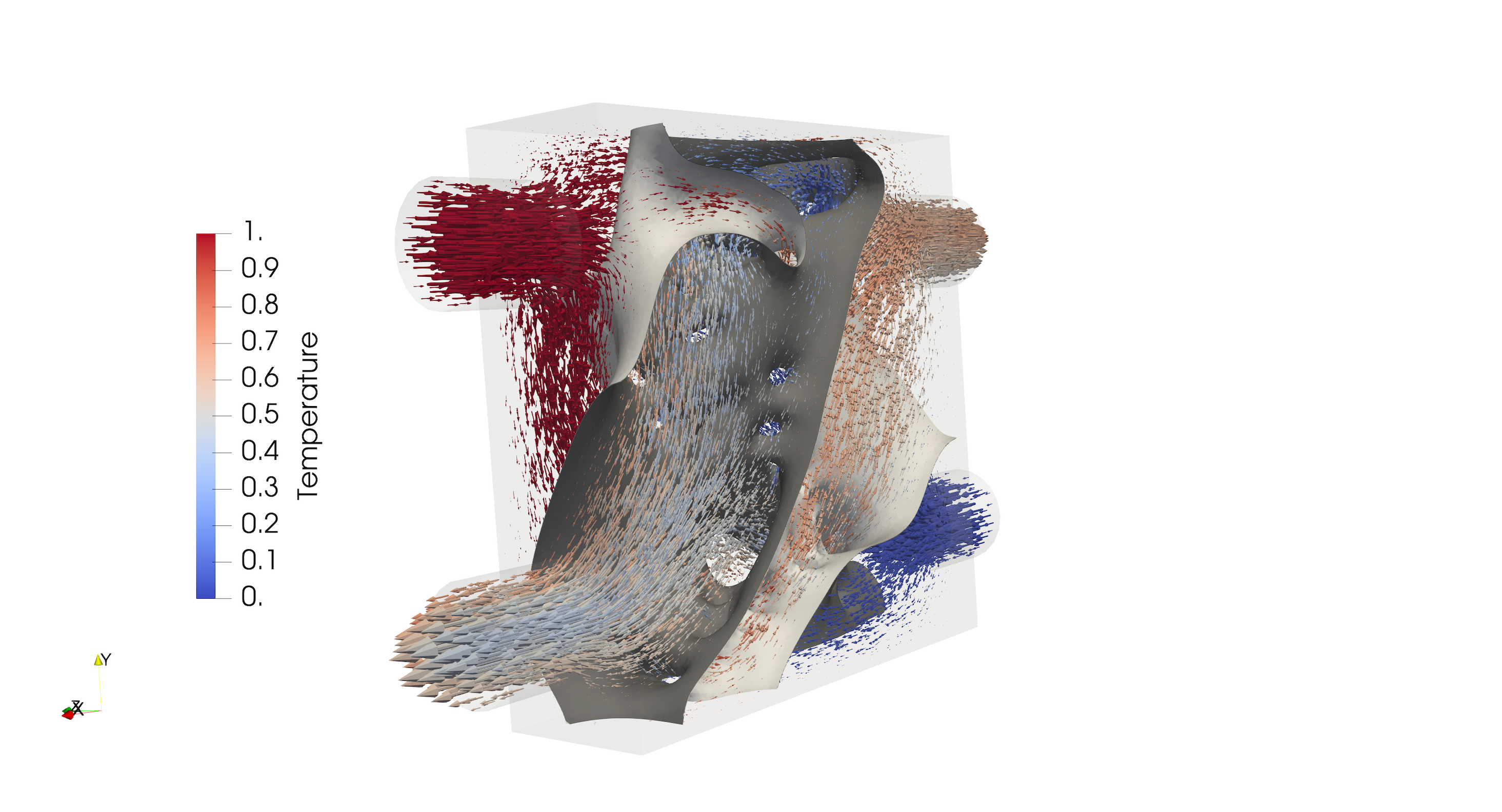}
	\caption{Optimized design, flow field and temperature for the counter flow heat exchanger with $P_{\text{drop}} = 2.0$ and $Re=10$, back view.
		Cold membrane side in black, hot side in white.}
	\label{fig:counter_diffuse}
\end{figure}
\begin{figure}[h!]
	\centering
	\includegraphics[scale=\imagesize]{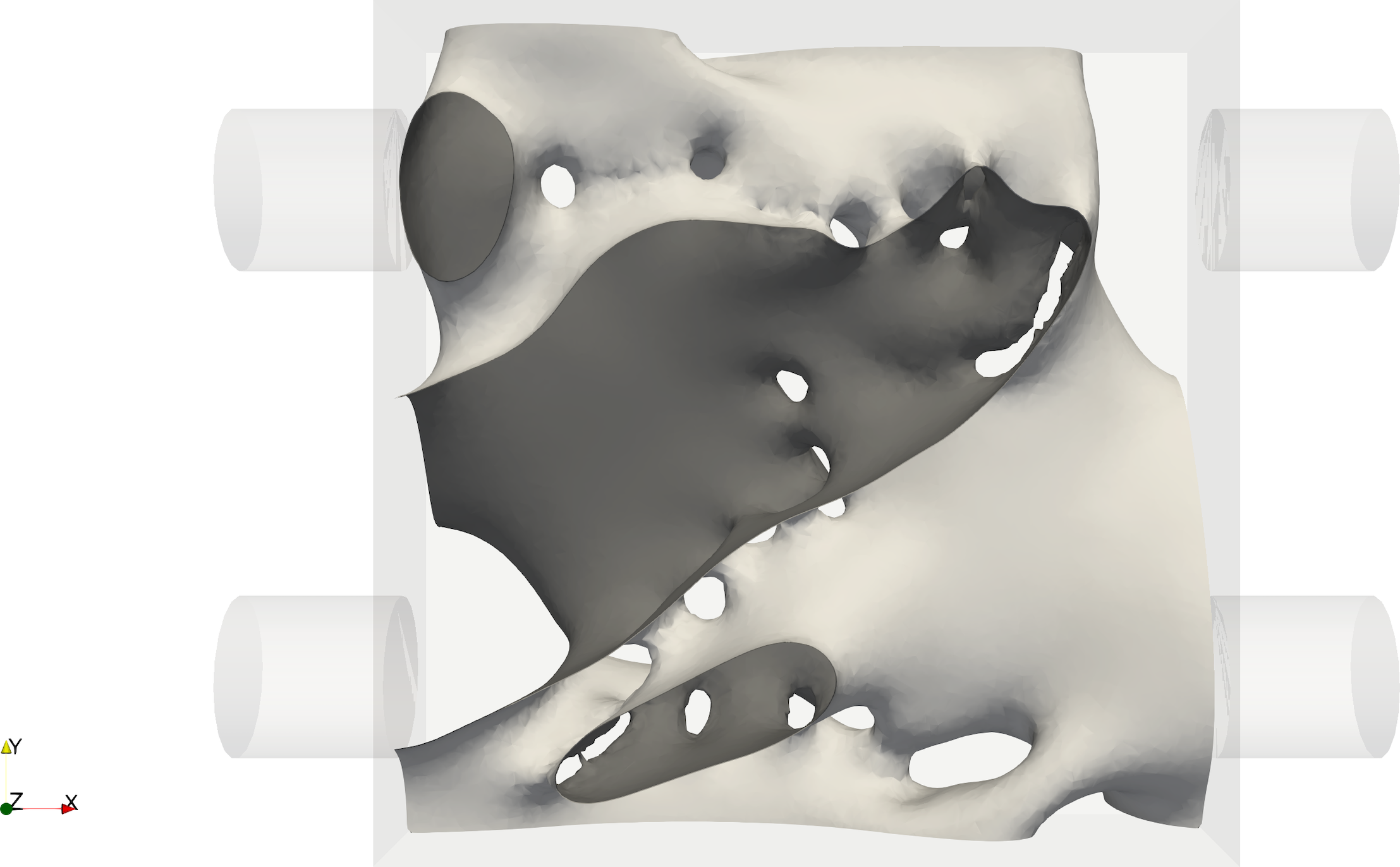}
	\caption{Optimized design for the counter flow heat exchanger with $P_{\text{drop}} = 2.0$ and $Re=10$.
		Cold membrane side in black, hot side in white.}
	\label{fig:counter_side}
\end{figure}

\begin{figure}[h!]
	\centering
	\includegraphics[scale=\imagesize]{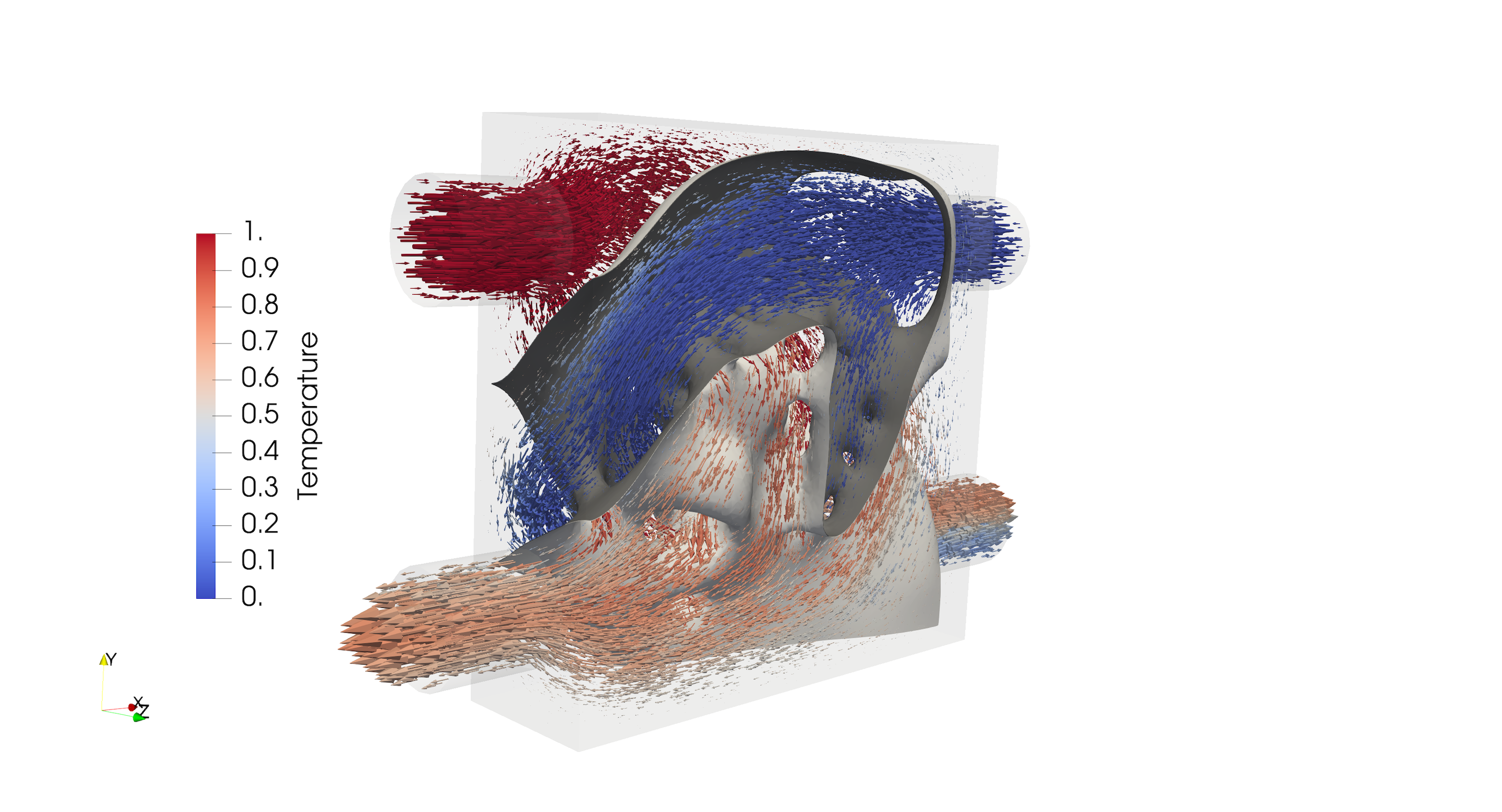}
	\caption{Optimized design, flow field and temperature for the U-flow heat exchanger with $P_{\text{drop}} = 2.0$ and $Re=10$, front view.
		Cold membrane side in black, hot side in white.}
	\label{fig:u_flow_clear}
\end{figure}
\begin{figure}[h!]
	\centering
	\includegraphics[scale=\imagesize]{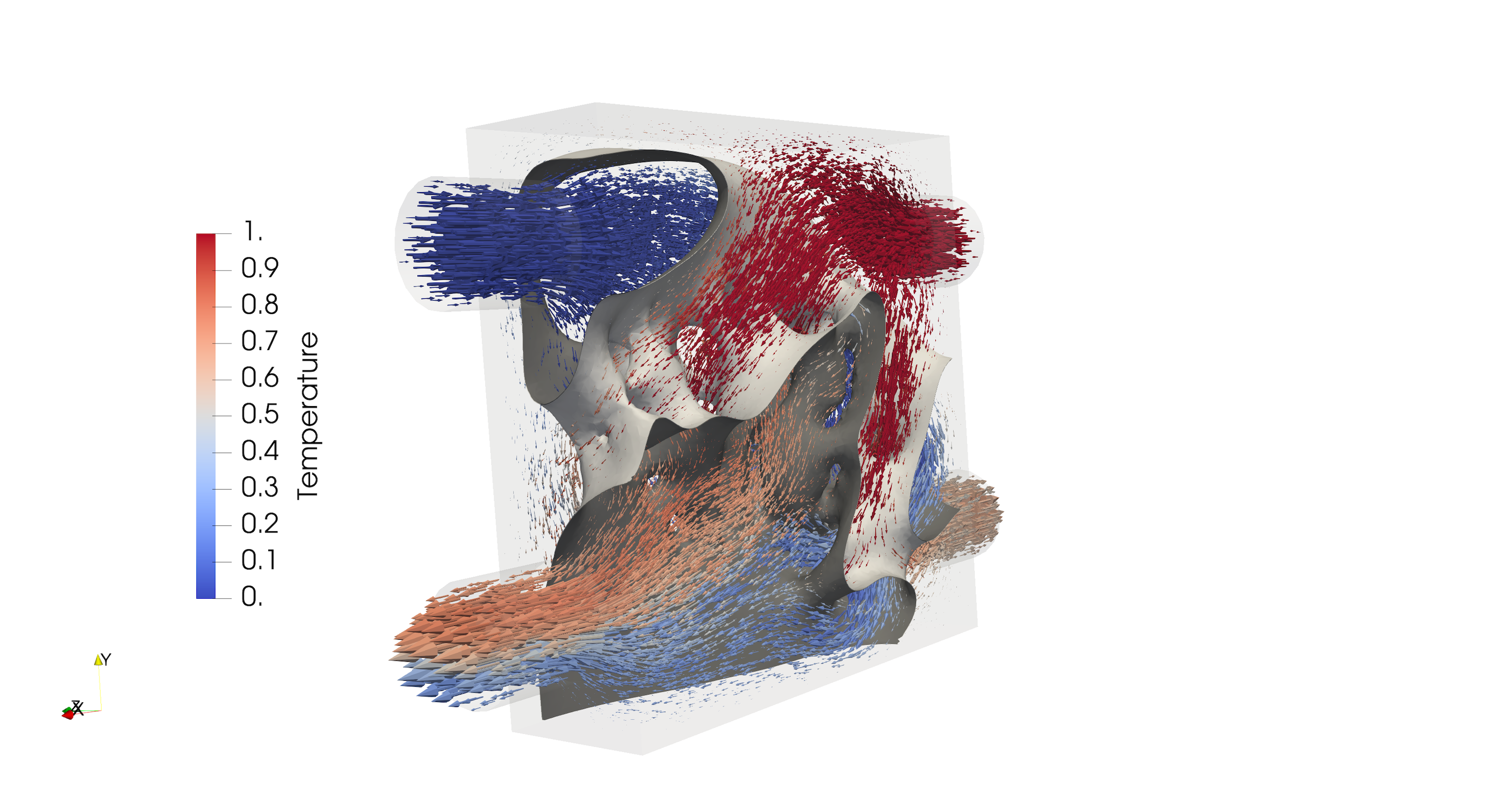}
	\caption{Optimized designs, flow field and temperature for the U-flow heat exchanger with $P_{\text{drop}} = 2.0$ and $Re=10$, back view.
		Cold membrane side in black, hot side in white.}
	\label{fig:u_flow_diffuse}
\end{figure}
\begin{figure}[h!]
	\centering
	\includegraphics[scale=\imagesize]{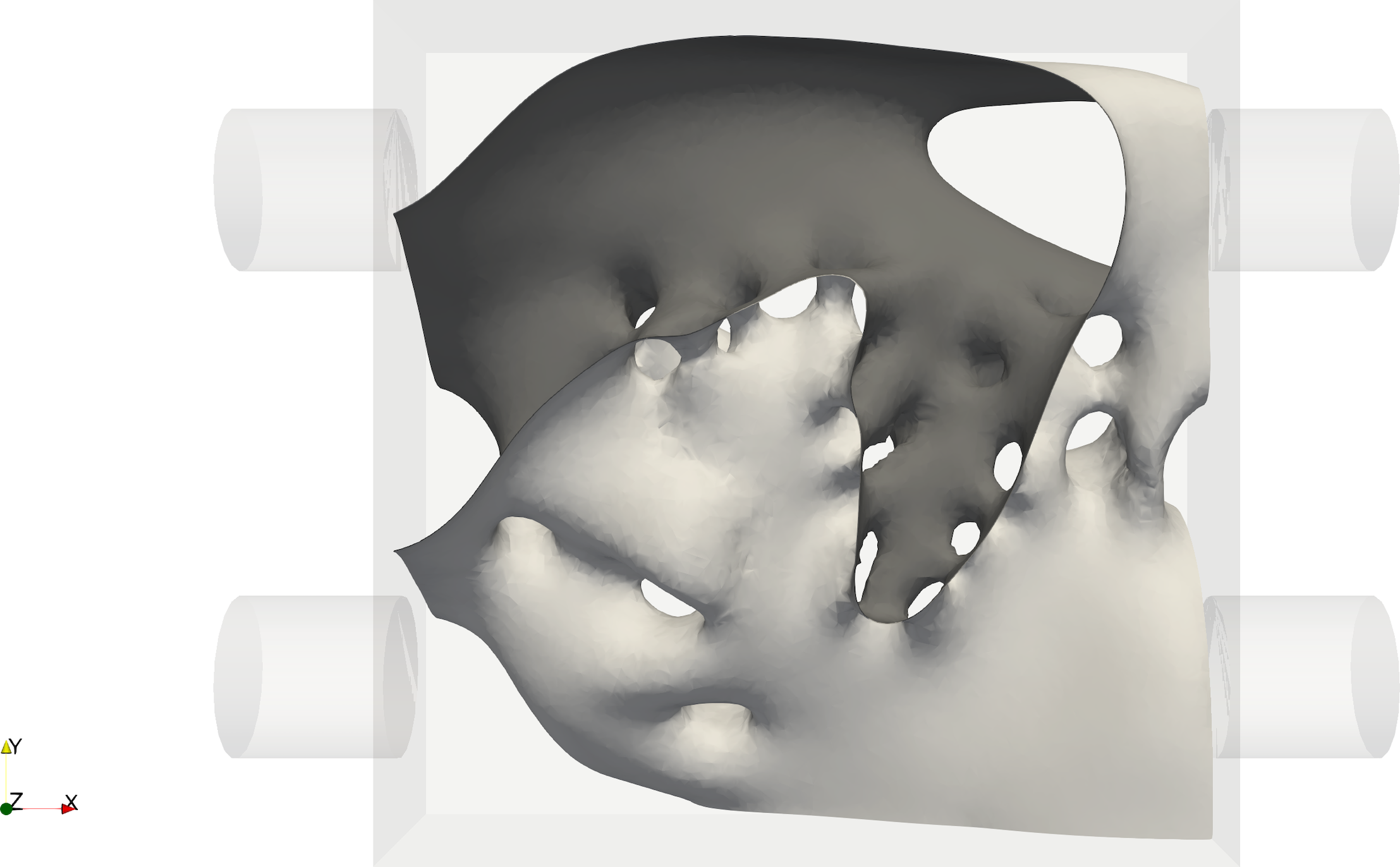}
	\caption{Optimized design for the U-flow heat exchanger with $P_{\text{drop}} = 2.0$ and $Re=10$.
		Cold membrane side in black, hot side in white.}
	\label{fig:u_flow_side}
\end{figure}

\begin{figure}[h!]
	\centering
	\includegraphics[scale=\imagesize]{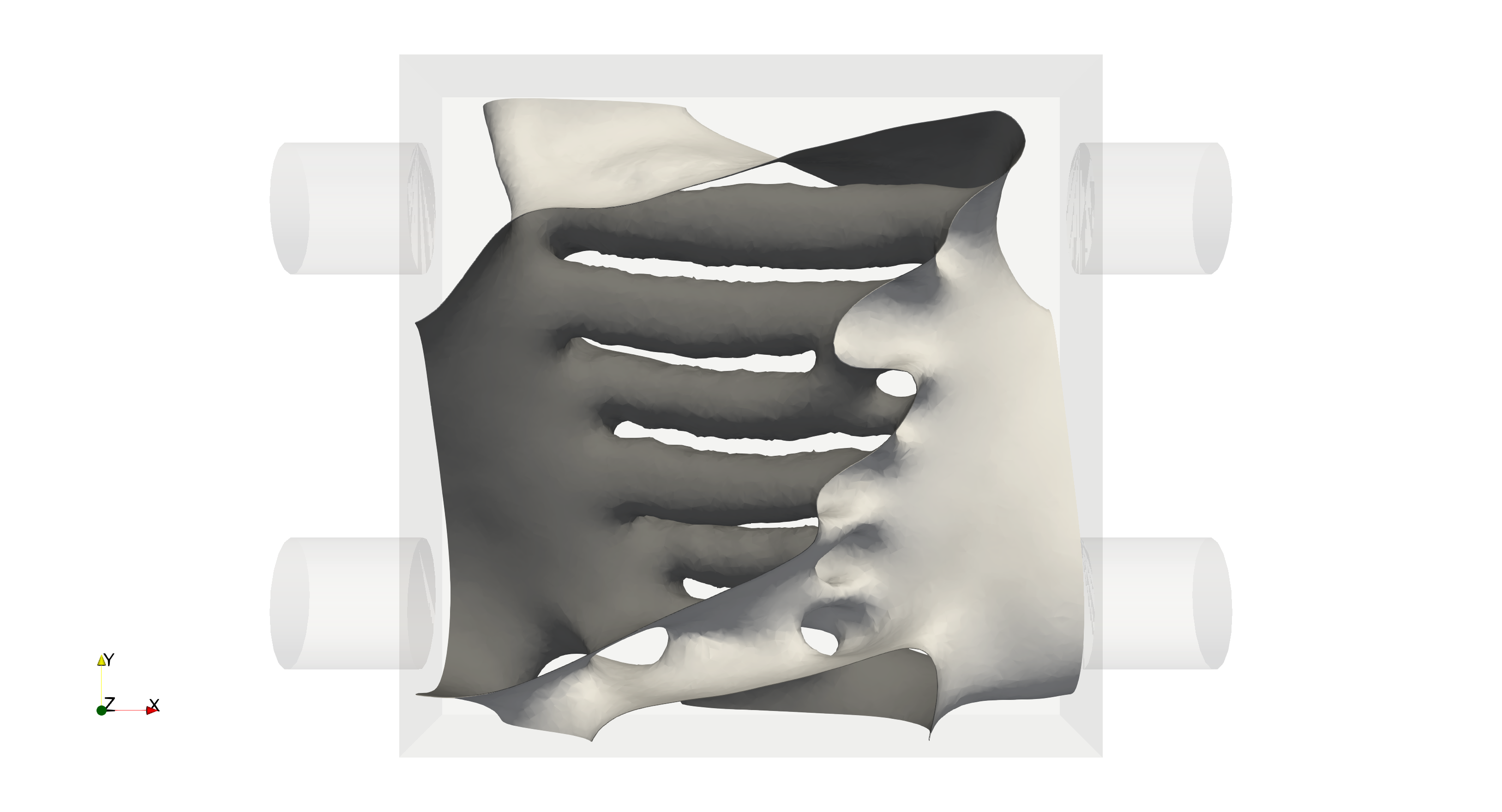}
	\caption{Optimized design for the parallel heat exchanger with $P_{\text{drop}} = 1.5$ and $Re=10$.
		Cold membrane side in black, hot side in white.}
	\label{fig:parallel_side_15}
\end{figure}

The next set of examples enforce different pressure drop constraint values, i.e. $P_{\text{drop}}=2.0, 1.5~\text{and}~1.0$.
Greater pressure drop allows for more narrow channels so as to increase the interface channel area between the fluids and hence maximize heat transfer,
cf. Table \ref{tab:cost_function_values_pressure}.
Notably, higher pressure drop allows the formation of connection channels perpendicular to the main flow, cf. Figures \ref{fig:parallel_side} and \ref{fig:parallel_side_15}

\begin{figure}[h!]
	\centering
	\includegraphics[scale=\imagesize]{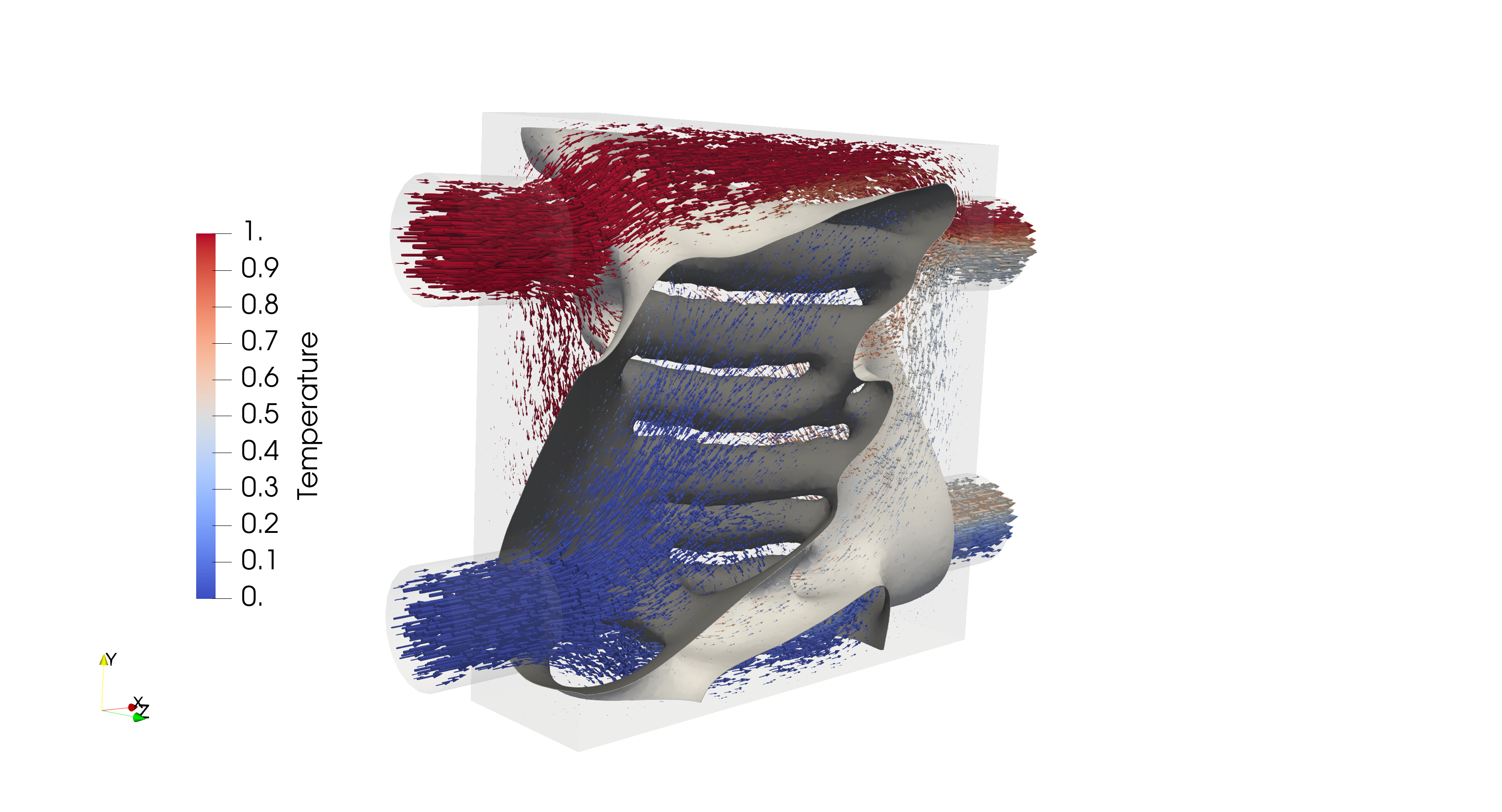}
	\caption{Optimized design, flow field and temperature for the parallel heat exchanger with $P_{\text{drop}} = 1.5$ and $Re=10$, front view.
		Cold membrane side in black, hot side in white.}
	\label{fig:parallel_clear_15}
\end{figure}
\begin{figure}[h!]
	\centering
	\includegraphics[scale=\imagesize]{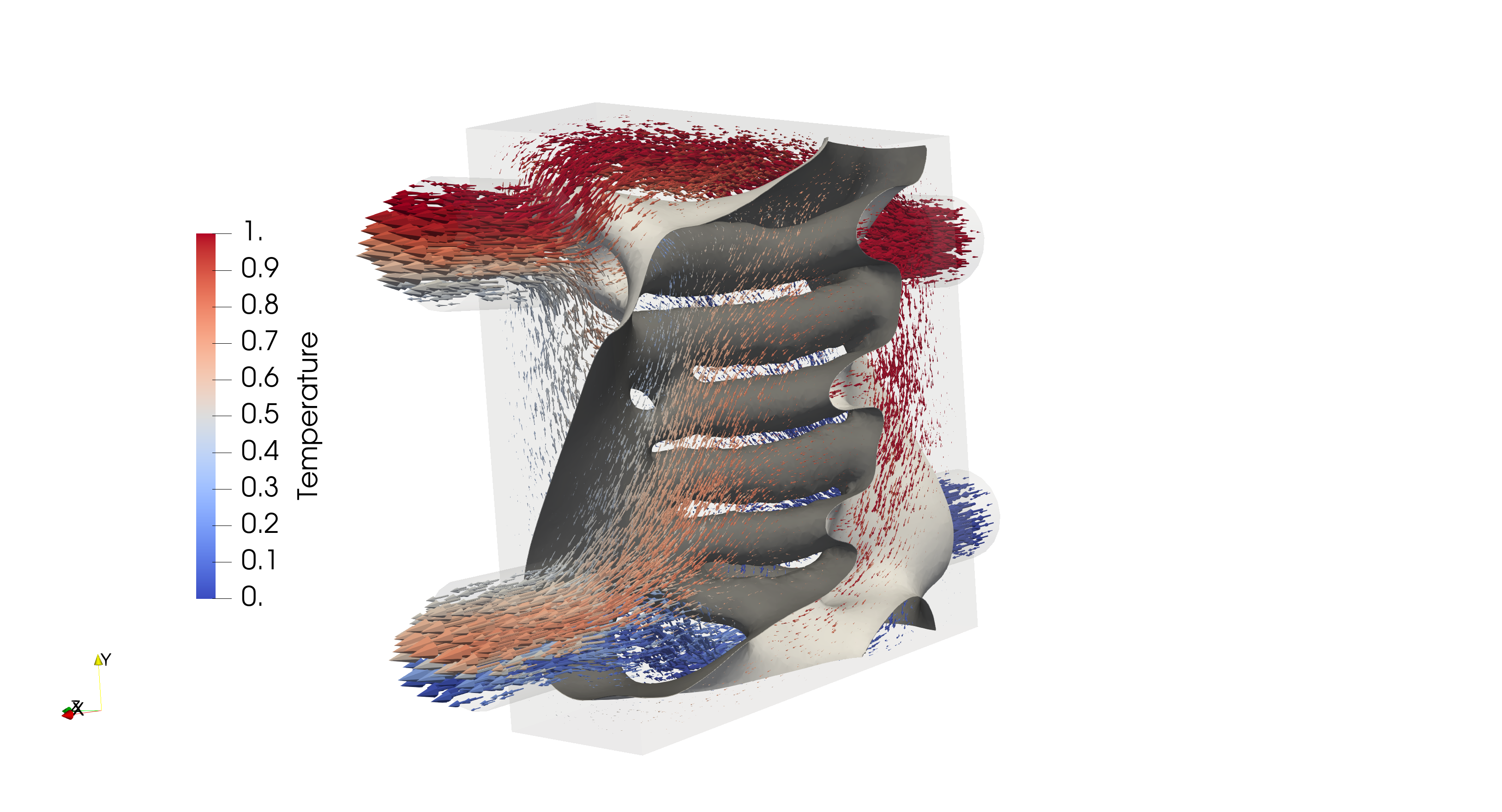}
	\caption{Optimized design, flow field and temperature for the parallel heat exchanger with $P_{\text{drop}} = 1.5$ and $Re=10$, back view.
		Cold membrane side in black, hot side in white.}
	\label{fig:parallel_diffuse_15}
\end{figure}
\begin{figure}[h!]
	\centering
	\includegraphics[scale=\imagesize]{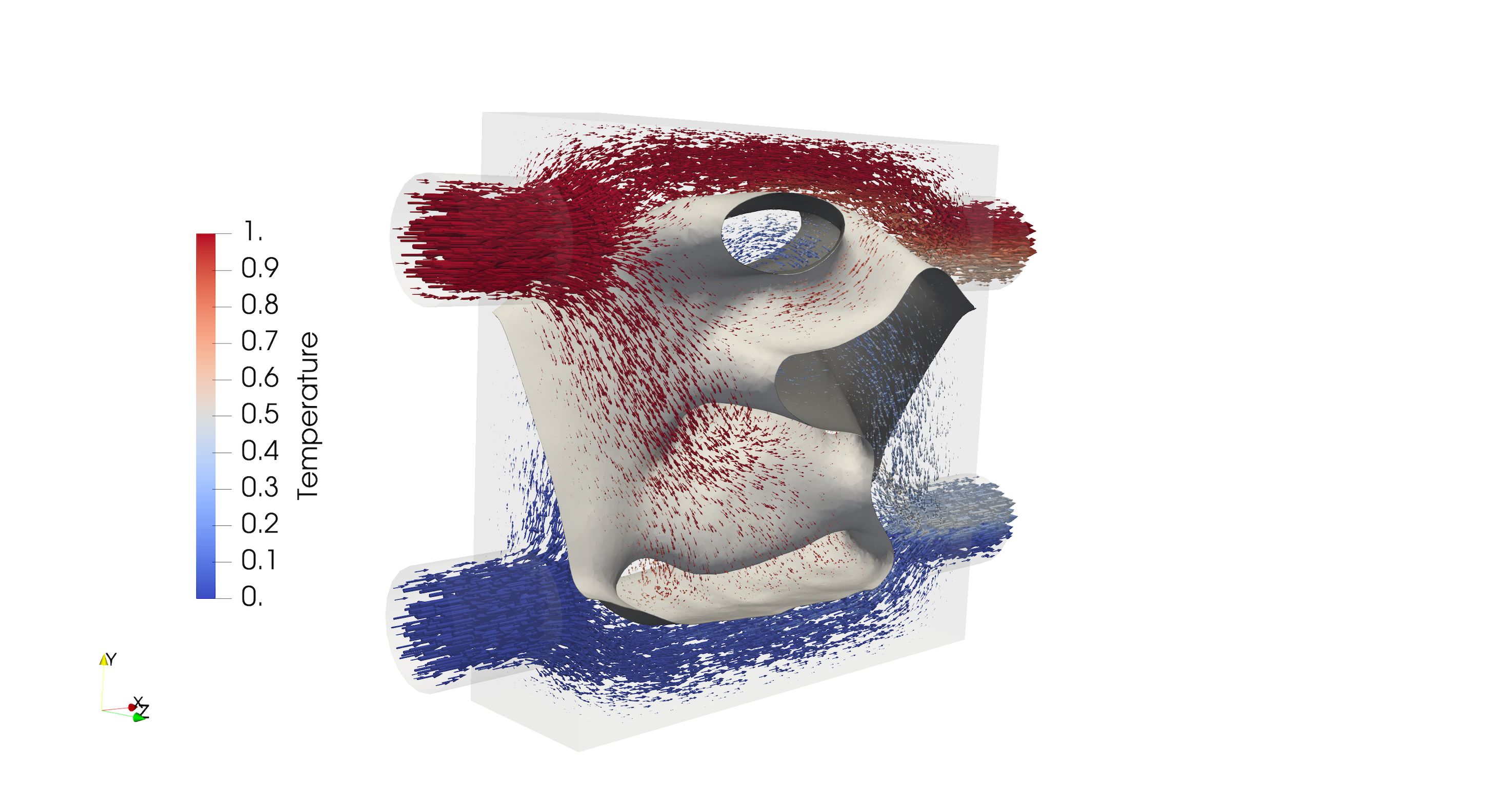}
	\caption{Optimized design, flow field and temperature for the parallel heat exchanger with $P_{\text{drop}} = 1.0$ and $Re=10$, front view.
		Cold membrane side in black, hot side in white.}
	\label{fig:parallel_clear_10}
\end{figure}
\begin{figure}[h!]
	\centering
	\includegraphics[scale=\imagesize]{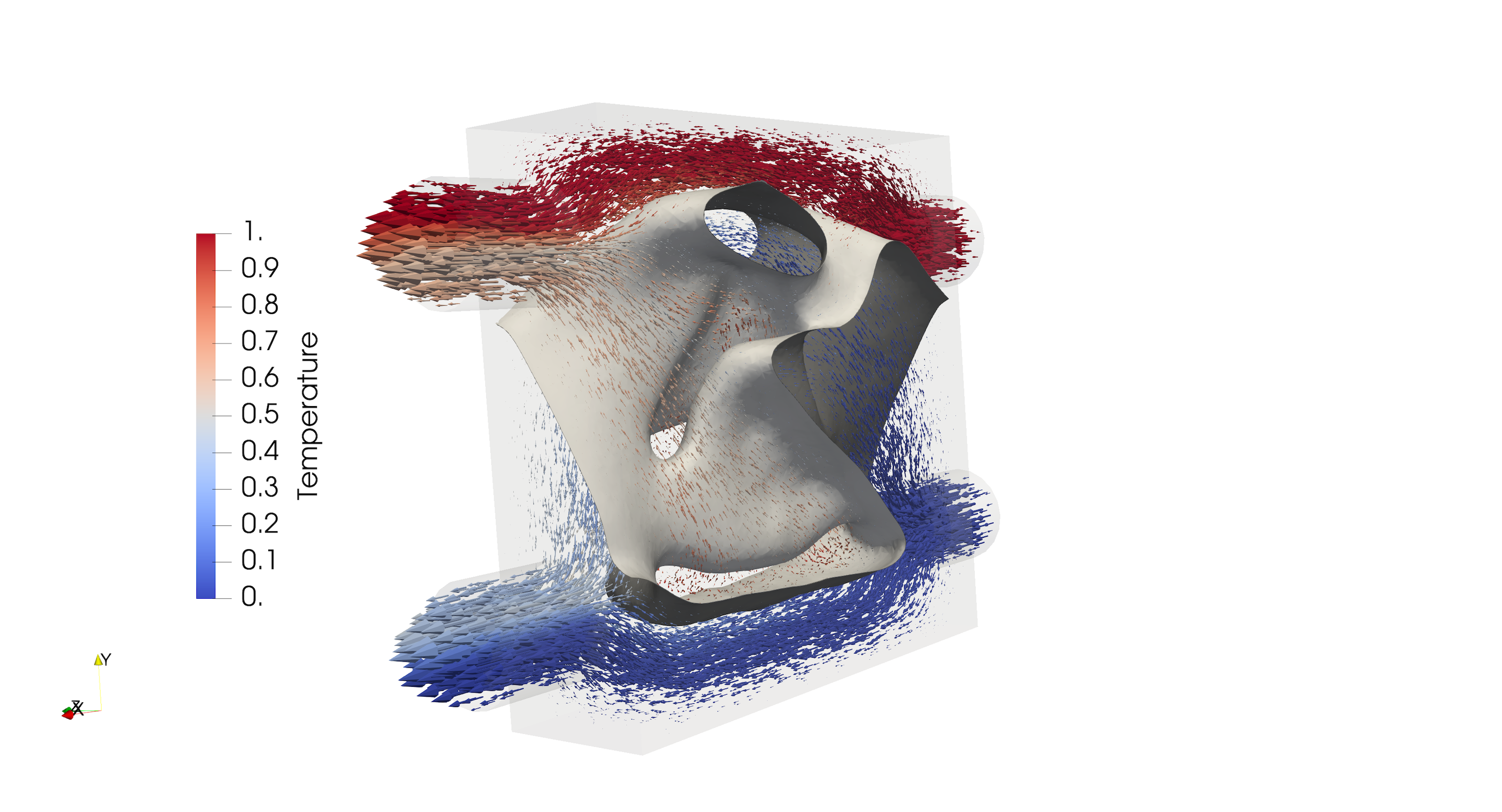}
	\caption{Optimized design, flow field and temperature for the parallel heat exchanger with $P_{\text{drop}} = 1.0$ and $Re=10$, back view.
		Cold membrane side in black, hot side in white.}
	\label{fig:parallel_diffuse_10}
\end{figure}
\begin{figure}[h!]
	\centering
	\includegraphics[scale=\imagesize]{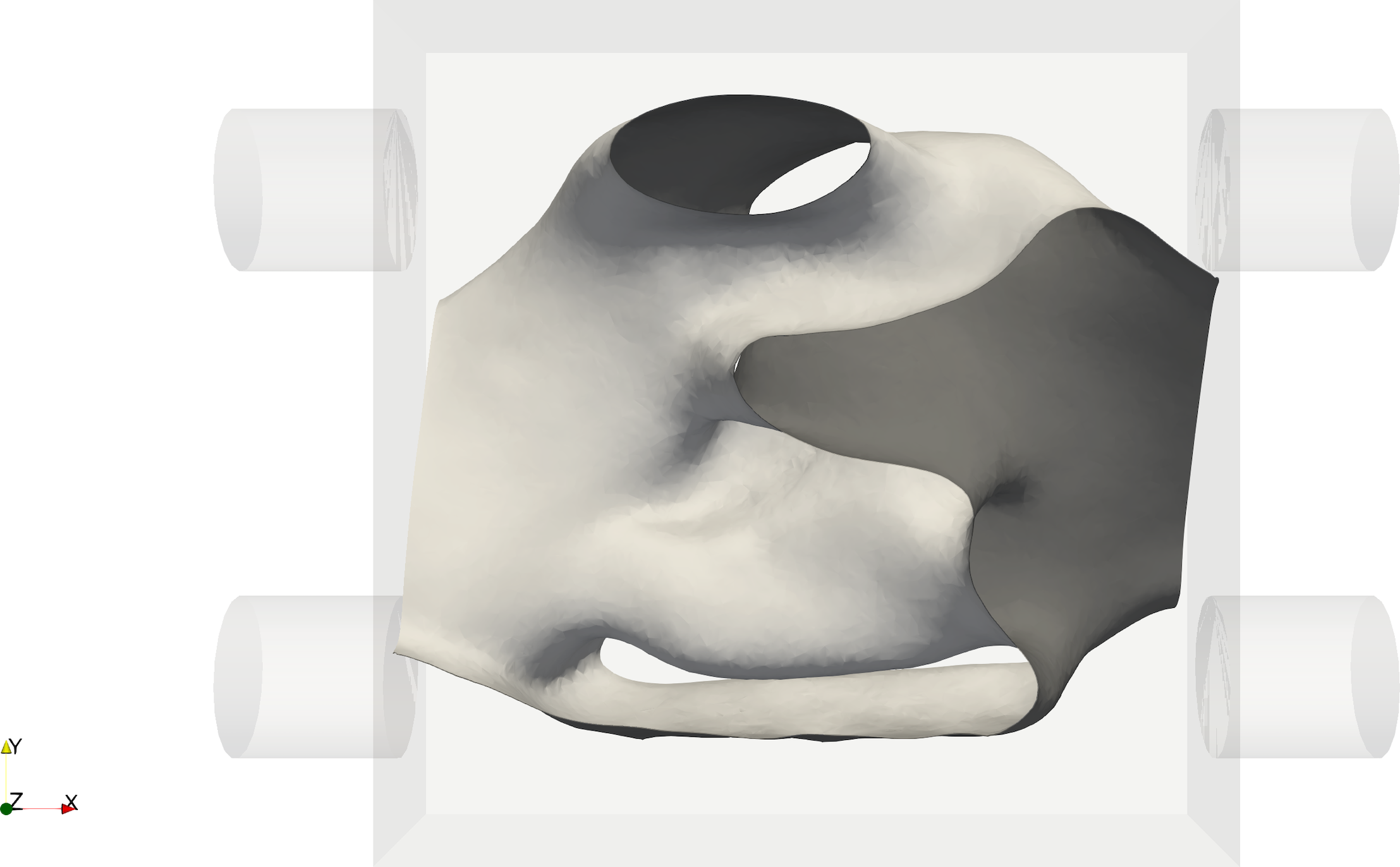}
	\caption{Optimized design for the parallel heat exchanger with $P_{\text{drop}} = 1.0$ and $Re=10$.
		Cold membrane side in black, hot side in white.}
	\label{fig:parallel_side_10}
\end{figure}
\begin{table}[]
	\centering
	\begin{tabular}{l|l|l|l}
		$P_{\text{drop}}$ & 2.0    & 1.5    & 1.0    \\
		\hline
		Cost function     & -0.489 & -0.337 & -0.211 \\
	\end{tabular}
	\caption{Cost function values for the optimized designs for $P_{\text{drop}}=2.0, 1.5~\text{and}~1.0$ }
	\label{tab:cost_function_values_pressure}
\end{table}

\begin{figure}[h!]
	\centering
	\includegraphics[scale=\imagesize]{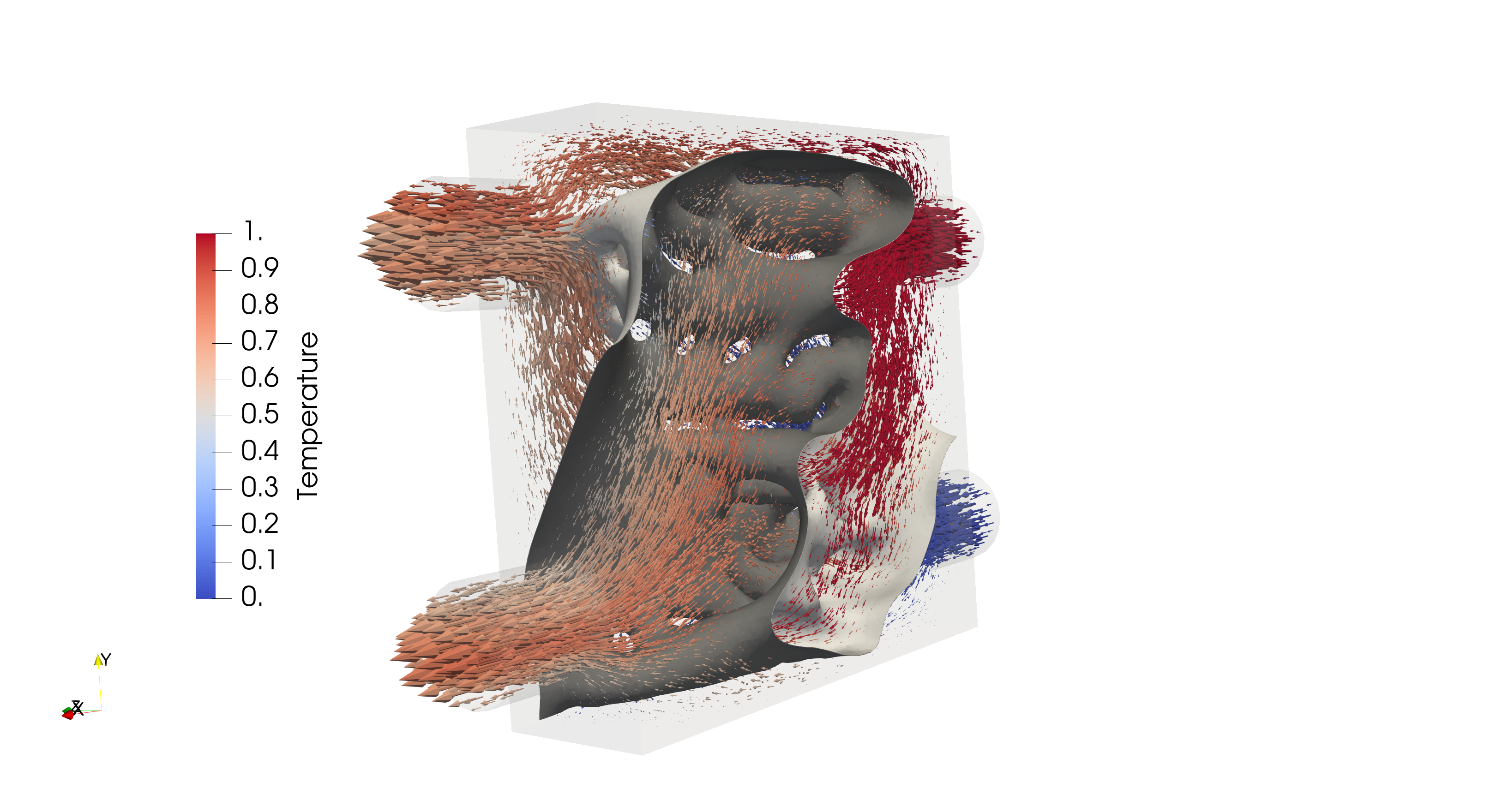}
	\caption{Optimized designs, flow field and temperature for the parallel heat exchanger with $P_{\text{drop}} = 1.0$ and $Re=20$.
		Cold membrane side in black, hot side in white.}
	\label{fig:parallel_diffuse_nu_005}
\end{figure}

\begin{figure}[h!]
	\centering
	\includegraphics[scale=\imagesize]{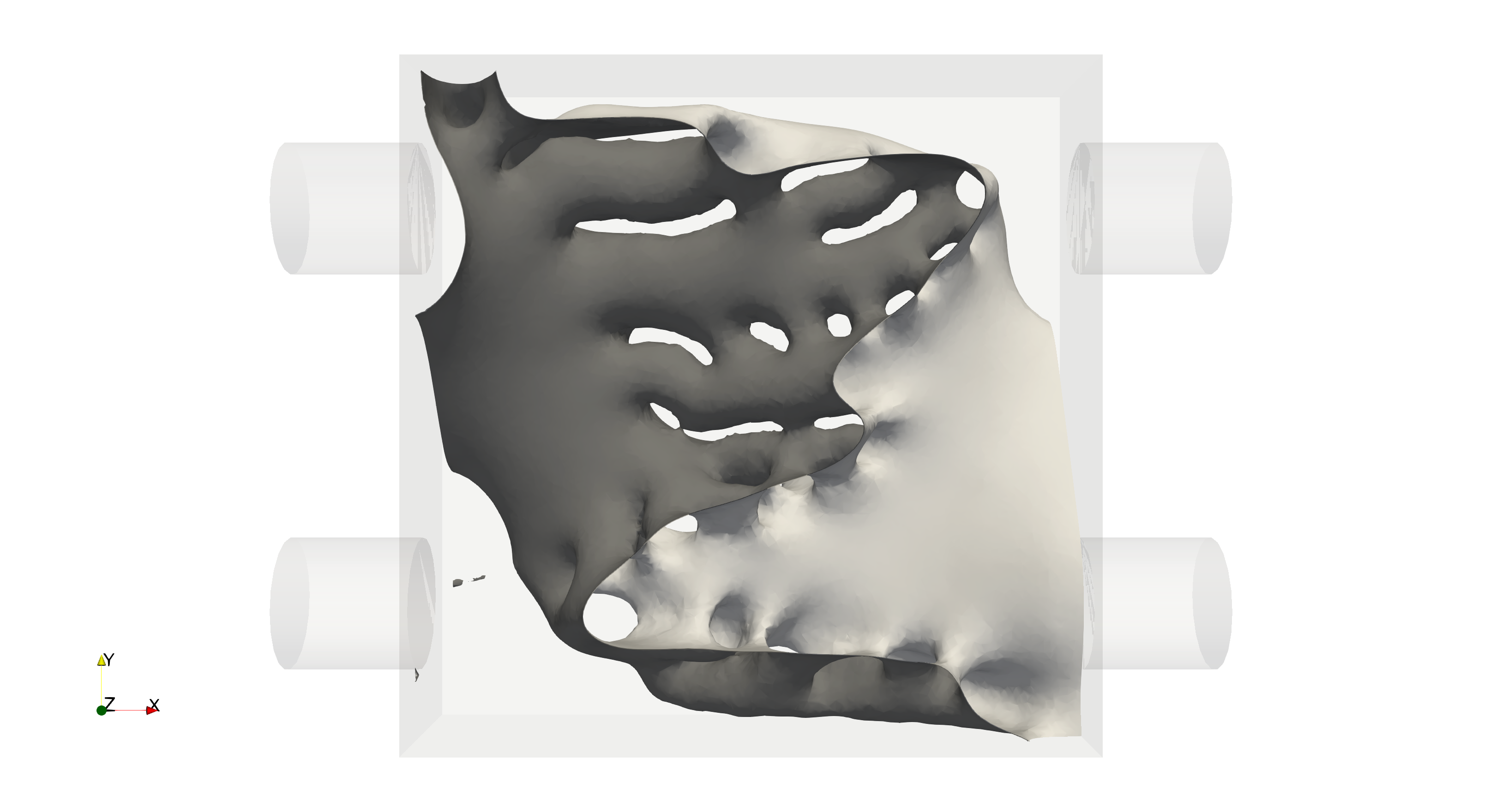}
	\caption{Optimized designs, flow field and temperature for the parallel heat exchanger with $P_{\text{drop}} = 2.0$ and $Re=20$.
		Cold membrane side in black, hot side in white.}
	\label{fig:parallel_side_nu_005}
\end{figure}

Next we increase the Reynolds number to $Re=20$ in a parallel configuration with $P_{\text{drop}} = 1.0$.
The optimized design, cf. Figures \ref{fig:parallel_diffuse_nu_005}-\ref{fig:parallel_side_nu_005} exhibits the complexity of the $Re=10$ and $P_{\text{drop}}=2.0$ design.
Indeed, a higher Reynolds number translates into a lower pressure drop across the heat exchanger, giving the optimizer the freedom to increase the complexity.
Unfortunately, examples with still higher Reynolds number require better preconditioners and parameter tuning of the flow solvers which we will address in our future work.
Table \ref{tab:cost_function_values_reynolds} denotes the cost function values for the $Re=10$ and $Re=20$ designs with $P_{\text{drop}} = 1.0$

\begin{figure}[h!]
	\centering
	\includegraphics[scale=\imagesize]{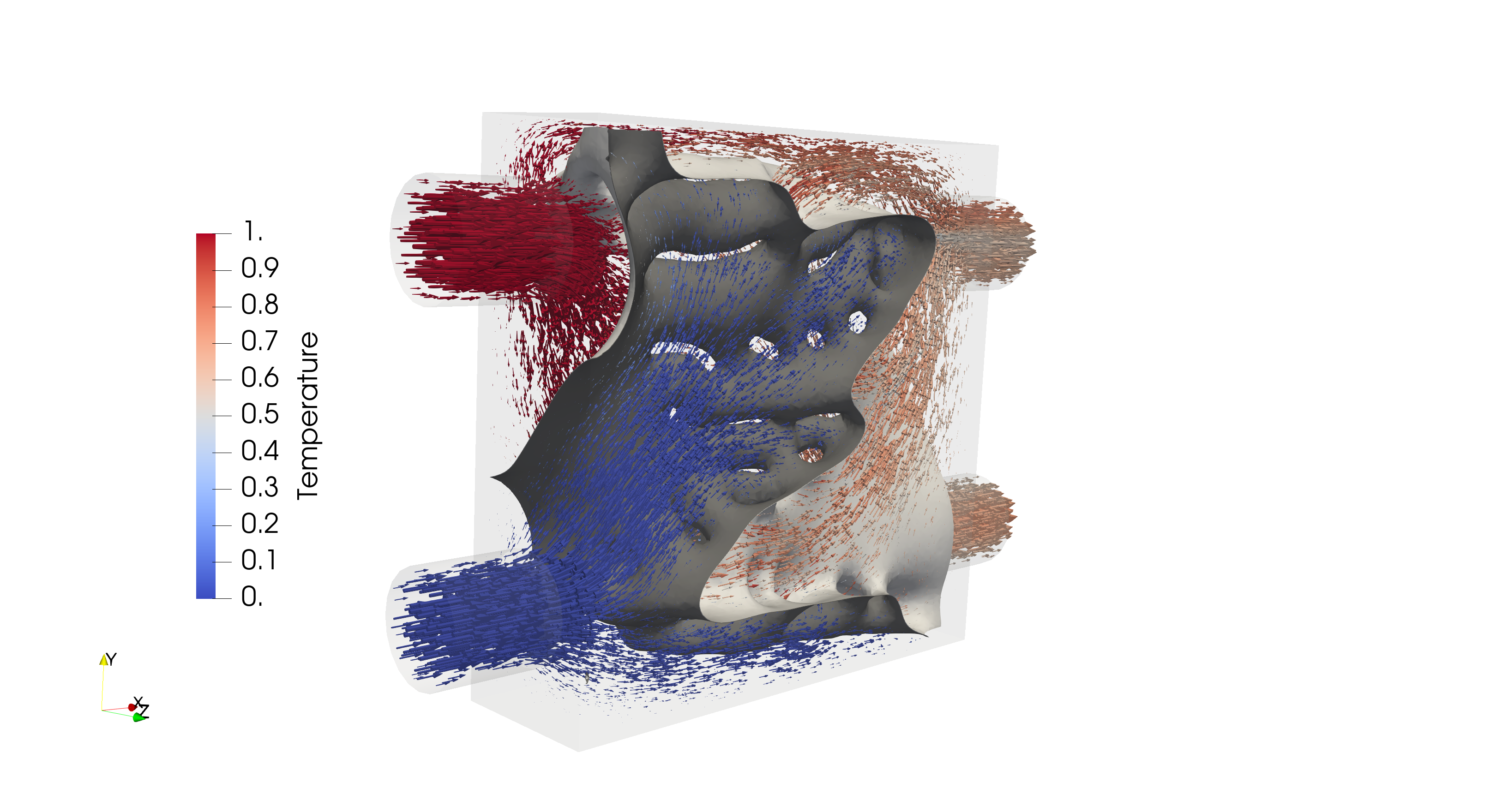}
	\caption{Optimized designs, flow field and temperature for the parallel heat exchanger with $P_{\text{drop}} = 1.0$ and $Re=20$.
		Cold membrane side in black, hot side in white.}
	\label{fig:parallel_clear_nu_005}
\end{figure}
\begin{table}[]
	\centering
	\begin{tabular}{l|l|l|l}
		Reynolds number & 10.0   & 20.0   \\
		\hline
		Cost function   & -0.211 & -0.498 \\
	\end{tabular}
	\caption{Cost function values for the optimized designs for $Re=10~\text{and}~20$ }
	\label{tab:cost_function_values_reynolds}
\end{table}

\section{Conclusions}
\label{sec:conclusion}
We design heat exchangers using the level-set method with an Ersatz approach to permit complete evolution of the topology.
We model the fluid flow over the entire domain and prevent the fluids from mixing by solving the Navier-Stokes equation for each fluid with a Brinkmann penalization to model the complementary fluid domain as solid.
The velocity fields drive the heat transport via a coupled advection-diffusion equation for the temperature.
Using the computed response, we proceed to maximize the heat transfer from the hot to cold fluid with a constraint on the pressure drop across each fluid channel.
We solve the optimization problem for several heat exchanger configurations, different pressure constraints and different flow conditions.
Our designs increase in complexity as we relax the pressure drop constraint or increase the Reynolds number.
Despite the symmetry of the boundary conditions and geometry, for the initial design used in this study, the optimizer takes advantage of the full domain to create asymmetric optimized designs.
Indeed, the optimizer creates two main channels that cross each other over the XY plane.
More studies with symmetric initial designs could confirm the optimality of the asymmetric designs.
We have also assumed an infinitesimally thin interface separate the fluids.
A channel wall with finite thickness and distinct conductivity properties will provide more predictive designs.
Lastly, research on a better preconditioner will allow us to obtain designs with higher Reynolds number.
Future work will also use immersed finite element methods such as the CutFEM \citep{burman2015cutfem} or XFEM \citep{belytschko2003structured} to model the interface with more precision.

\section*{Replication of results}
The software used is archived in \cite{zenodo/Firedrake-20211025.1}.
It additionally requires the library Lestofire \cite{miguel_salazar_de_troya_2021_5597909}.
We offer ourselves to help any user interested in replicating the results.

\section*{Conflict of interest.}
The authors declare that they have no conflict of interest.

\section*{Funding information}
This work was performed under the auspices of the  U.S.  Department  of  Energy  by  Lawrence  Livermore  National Laboratory under Contract DE-AC52-07NA27344 and was supported by the LLNL-LDRD program under project number 19-ERD-035. LLNL Release Number LLNL-JRNL-816310.
\bibliographystyle{elsarticle-num}
\bibliography{heat.bib}

\begin{thebibliography}{10}
\expandafter\ifx\csname url\endcsname\relax
  \def\url#1{\texttt{#1}}\fi
\expandafter\ifx\csname urlprefix\endcsname\relax\def\urlprefix{URL }\fi
\expandafter\ifx\csname href\endcsname\relax
  \def\href#1#2{#2} \def\path#1{#1}\fi

\bibitem{Patel2019}
V.~K. Patel, V.~J. Savsani, M.~A. Tawhid,
  \href{https://doi.org/10.1007/978-3-030-10477-1_3}{Thermal Design and
  Optimization of Heat Exchangers}, Springer International Publishing, Cham,
  2019, pp. 33--98.
\newblock \href {https://doi.org/10.1007/978-3-030-10477-1_3}
  {\path{doi:10.1007/978-3-030-10477-1_3}}.
\newline\urlprefix\url{https://doi.org/10.1007/978-3-030-10477-1_3}

\bibitem{dixit2015review}
T.~Dixit, I.~Ghosh, Review of micro-and mini-channel heat sinks and heat
  exchangers for single phase fluids, Renewable and Sustainable Energy Reviews
  41 (2015) 1298--1311.

\bibitem{awais2018heat}
M.~Awais, A.~A. Bhuiyan, Heat and mass transfer for compact heat exchanger
  (chxs) design: A state-of-the-art review, International Journal of Heat and
  Mass Transfer 127 (2018) 359--380.

\bibitem{foli2006optimization}
K.~Foli, T.~Okabe, M.~Olhofer, Y.~Jin, B.~Sendhoff, Optimization of micro heat
  exchanger: Cfd, analytical approach and multi-objective evolutionary
  algorithms, International Journal of Heat and Mass Transfer 49~(5-6) (2006)
  1090--1099.

\bibitem{bendsoe}
M.~P. Bends{\o}e, Optimal shape design as a material distribution problem,
  Structural optimization 1~(4) (1989) 193--202.

\bibitem{wang2003level}
M.~Y. Wang, X.~Wang, D.~Guo, A level set method for structural topology
  optimization, Computer methods in applied mechanics and engineering 192~(1-2)
  (2003) 227--246.

\bibitem{allaire}
G.~Allaire, F.~Jouve, A.-M. Toader, A level-set method for shape optimization,
  Comptes Rendus Mathematique 334~(12) (2002) 1125--1130.

\bibitem{SETHIAN2000489}
J.~Sethian, A.~Wiegmann,
  \href{http://www.sciencedirect.com/science/article/pii/S0021999100965811}{Structural
  boundary design via level set and immersed interface methods}, Journal of
  Computational Physics 163~(2) (2000) 489 -- 528.
\newblock \href {https://doi.org/https://doi.org/10.1006/jcph.2000.6581}
  {\path{doi:https://doi.org/10.1006/jcph.2000.6581}}.
\newline\urlprefix\url{http://www.sciencedirect.com/science/article/pii/S0021999100965811}

\bibitem{borrvall}
T.~Borrvall, J.~Petersson, Topology optimization of fluids in stokes flow,
  International journal for numerical methods in fluids 41~(1) (2003) 77--107.

\bibitem{gersborg2005topology}
A.~Gersborg-Hansen, O.~Sigmund, R.~B. Haber, Topology optimization of channel
  flow problems, Structural and Multidisciplinary Optimization 30~(3) (2005)
  181--192.

\bibitem{kiziltas2003topology}
G.~Kiziltas, D.~Psychoudakis, J.~L. Volakis, N.~Kikuchi, Topology design
  optimization of dielectric substrates for bandwidth improvement of a patch
  antenna, IEEE Transactions on Antennas and Propagation 51~(10) (2003)
  2732--2743.

\bibitem{zhou2010level}
S.~Zhou, W.~Li, Q.~Li, Level-set based topology optimization for
  electromagnetic dipole antenna design, Journal of Computational Physics
  229~(19) (2010) 6915--6930.

\bibitem{wadbro}
E.~Wadbro, M.~Berggren, Topology optimization of an acoustic horn, Computer
  Methods in Applied Mechanics and Engineering 196 (2006) 420–436.
\newblock \href {https://doi.org/10.1016/j.cma.2006.05.005}
  {\path{doi:10.1016/j.cma.2006.05.005}}.

\bibitem{deaton2014survey}
J.~D. Deaton, R.~V. Grandhi, A survey of structural and multidisciplinary
  continuum topology optimization: post 2000, Structural and Multidisciplinary
  Optimization 49~(1) (2014) 1--38.

\bibitem{van2013level}
N.~P. van Dijk, K.~Maute, M.~Langelaar, F.~Van~Keulen, Level-set methods for
  structural topology optimization: a review, Structural and Multidisciplinary
  Optimization 48~(3) (2013) 437--472.

\bibitem{alexandersen2020review}
J.~Alexandersen, C.~S. Andreasen, A review of topology optimisation for
  fluid-based problems, Fluids 5~(1) (2020) 29.

\bibitem{gersborg2006topology}
A.~Gersborg-Hansen, M.~P. Bends{\o}e, O.~Sigmund, Topology optimization of heat
  conduction problems using the finite volume method, Structural and
  multidisciplinary optimization 31~(4) (2006) 251--259.

\bibitem{li1999shape}
Q.~Li, G.~P. Steven, O.~M. Querin, Y.~Xie, Shape and topology design for heat
  conduction by evolutionary structural optimization, International Journal of
  Heat and Mass Transfer 42~(17) (1999) 3361--3371.

\bibitem{dede2009multiphysics}
E.~M. Dede, Multiphysics topology optimization of heat transfer and fluid flow
  systems, in: In Proceedings of the COMSOL Conference 2009, , Boston, MA, USA,
  8–10 October 2009., 2009.

\bibitem{yoon2010topological}
G.~H. Yoon, Topological design of heat dissipating structure with forced
  convective heat transfer, Journal of Mechanical Science and Technology 24~(6)
  (2010) 1225--1233.

\bibitem{alexandersen2014topology}
J.~Alexandersen, N.~Aage, C.~S. Andreasen, O.~Sigmund, Topology optimisation
  for natural convection problems, International Journal for Numerical Methods
  in Fluids 76~(10) (2014) 699--721.

\bibitem{heatreview}
T.~Dbouk, A review about the engineering design of optimal heat transfer
  systems using topology optimization, Applied Thermal Engineering 112 (2017)
  841--854.

\bibitem{haertel2017fully}
J.~H. Haertel, G.~F. Nellis, A fully developed flow thermofluid model for
  topology optimization of 3d-printed air-cooled heat exchangers, Applied
  thermal engineering 119 (2017) 10--24.

\bibitem{kobayashi2019}
H.~Kobayashi, K.~Yaji, S.~Yamasaki, K.~Fujita, Freeform winglet design of
  fin-and-tube heat exchangers guided by topology optimization, Applied Thermal
  Engineering 161 (2019) 114020.
\newblock \href {https://doi.org/10.1016/j.applthermaleng.2019.114020}
  {\path{doi:10.1016/j.applthermaleng.2019.114020}}.

\bibitem{panas}
P.~Papazoglou, Topology optimization of heat exchangers (2015).

\bibitem{saviers2019design}
K.~R. Saviers, R.~Ranjan, R.~Mahmoudi, Design and validation of topology
  optimized heat exchangers, in: AIAA Scitech 2019 Forum, 2019, p. 1465.

\bibitem{Kobayashi2020}
H.~Kobayashi, K.~Yaji, S.~Yamasaki, K.~Fujita,
  \href{https://doi.org/10.1007/s00158-020-02736-8}{Topology design of
  two-fluid heat exchange}, Structural and Multidisciplinary Optimization (Oct
  2020).
\newblock \href {https://doi.org/10.1007/s00158-020-02736-8}
  {\path{doi:10.1007/s00158-020-02736-8}}.
\newline\urlprefix\url{https://doi.org/10.1007/s00158-020-02736-8}

\bibitem{Hoghoj}
L.~C. Høghøj, D.~R. Nørhave, J.~Alexandersen, O.~Sigmund, C.~S. Andreasen,
  \href{http://www.sciencedirect.com/science/article/pii/S0017931020334797}{Topology
  optimization of two fluid heat exchangers}, International Journal of Heat and
  Mass Transfer 163 (2020) 120543.
\newblock \href
  {https://doi.org/https://doi.org/10.1016/j.ijheatmasstransfer.2020.120543}
  {\path{doi:https://doi.org/10.1016/j.ijheatmasstransfer.2020.120543}}.
\newline\urlprefix\url{http://www.sciencedirect.com/science/article/pii/S0017931020334797}

\bibitem{feppon2}
F.~Feppon, G.~Allaire, C.~Dapogny, P.~Jolivet,
  \href{https://www.sciencedirect.com/science/article/pii/S0045782520308239}{Body-fitted
  topology optimization of {{2D}} and {{3D}} fluid-to-fluid heat exchangers},
  Computer Methods in Applied Mechanics and Engineering 376 (2021) 113638.
\newblock \href {https://doi.org/10.1016/j.cma.2020.113638}
  {\path{doi:10.1016/j.cma.2020.113638}}.
\newline\urlprefix\url{https://www.sciencedirect.com/science/article/pii/S0045782520308239}

\bibitem{Nunes:2020fg}
S.~P. Nunes, P.~Z. Culfaz-Emecen, G.~Z. Ramon, T.~Visser, G.~H. Koops, W.~Jin,
  M.~Ulbricht, {Thinking the future of membranes: Perspectives for advanced and
  new membrane materials and manufacturing processes}, Journal of Membrane
  Science 598 (2020) 117761.

\bibitem{Wang:2020jk}
Y.~Wang, D.~F. Ruiz~Diaz, K.~S. Chen, Z.~Wang, X.~C. Adroher, {Materials,
  technological status, and fundamentals of PEM fuel cells {\textendash} A
  review}, Materials Today 32 (2020) 178--203.

\bibitem{Park:2016bt}
M.~Park, J.~Ryu, W.~Wang, J.~Cho, {Material design and engineering of
  next-generation flow-battery technologies}, Nat. Rev. Mater. 2 (2017) 16080.

\bibitem{troyaTwoDimensionalTopology2021}
M.~Salazar~de Troya, D.~Tortorelli, V.~Beck,
  \href{https://www.scipedia.com/public/Troya_et_al_2021a}{Two {{Dimensional
  Topology Optimization}} of {{Heat Exchangers}} with the {{Density}} and
  {{Level}}-{{Set Methods}}}, in: 14th {{WCCM}}-{{ECCOMAS Congress}}, {CIMNE},
  2021.
\newblock \href {https://doi.org/10.23967/wccm-eccomas.2020.345}
  {\path{doi:10.23967/wccm-eccomas.2020.345}}.
\newline\urlprefix\url{https://www.scipedia.com/public/Troya_et_al_2021a}

\bibitem{WANG2018553}
Y.~Wang, Z.~Kang,
  \href{https://www.sciencedirect.com/science/article/pii/S0045782516319314}{A
  level set method for shape and topology optimization of coated structures},
  Computer Methods in Applied Mechanics and Engineering 329 (2018) 553--574.
\newblock \href {https://doi.org/10.1016/j.cma.2017.09.017}
  {\path{doi:10.1016/j.cma.2017.09.017}}.
\newline\urlprefix\url{https://www.sciencedirect.com/science/article/pii/S0045782516319314}

\bibitem{kreissl2011explicit}
S.~Kreissl, G.~Pingen, K.~Maute, An explicit level set approach for generalized
  shape optimization of fluids with the lattice boltzmann method, International
  Journal for Numerical Methods in Fluids 65~(5) (2011) 496--519.

\bibitem{ham2019automated}
D.~A. Ham, L.~Mitchell, A.~Paganini, F.~Wechsung, Automated shape
  differentiation in the unified form language, Structural and
  Multidisciplinary Optimization 60~(5) (2019) 1813--1820.

\bibitem{dokken2020automatic}
J.~S. Dokken, S.~K. Mitusch, S.~W. Funke, Automatic shape derivatives for
  transient pdes in fenics and firedrake, arXiv preprint arXiv:2001.10058
  (2020).

\bibitem{Mitusch2019}
S.~Mitusch, S.~Funke, J.~Dokken,
  \href{https://doi.org/10.21105/joss.01292}{dolfin-adjoint 2018.1: automated
  adjoints for fenics and firedrake}, Journal of Open Source Software 4~(38)
  (2019) 1292.
\newblock \href {https://doi.org/10.21105/joss.01292}
  {\path{doi:10.21105/joss.01292}}.
\newline\urlprefix\url{https://doi.org/10.21105/joss.01292}

\bibitem{laurain2018level}
A.~Laurain, A level set-based structural optimization code using fenics,
  Structural and Multidisciplinary Optimization 58~(3) (2018) 1311--1334.

\bibitem{delfour2011shapes}
M.~C. Delfour, J.-P. Zol{\'e}sio, Shapes and geometries: metrics, analysis,
  differential calculus, and optimization, Vol.~22, Siam, 2011.

\bibitem{choi1986design}
K.~K. Choi, H.~G. Seong, Design component method for sensitivity analysis of
  built-up structures, Journal of structural mechanics 14~(3) (1986) 379--399.

\bibitem{hiptmair2015comparison}
R.~Hiptmair, A.~Paganini, S.~Sargheini, Comparison of approximate shape
  gradients, BIT Numerical Mathematics 55~(2) (2015) 459--485.

\bibitem{berggren2010unified}
M.~Berggren, A unified discrete--continuous sensitivity analysis method for
  shape optimization, in: Applied and numerical partial differential equations,
  Springer, 2010, pp. 25--39.

\bibitem{gournay2006velocity}
F.~De~Gournay, Velocity extension for the level-set method and multiple
  eigenvalues in shape optimization, SIAM journal on control and optimization
  45~(1) (2006) 343--367.

\bibitem{burger2003framework}
M.~Burger, A framework for the construction of level set methods for shape
  optimization and reconstruction, Interfaces and Free boundaries 5~(3) (2003)
  301--329.

\bibitem{allaire2014multi}
G.~Allaire, C.~Dapogny, G.~Delgado, G.~Michailidis, Multi-phase structural
  optimization via a level set method, ESAIM: control, optimisation and
  calculus of variations 20~(2) (2014) 576--611.

\bibitem{dunning2015introducing}
P.~D. Dunning, H.~A. Kim, Introducing the sequential linear programming
  level-set method for topology optimization, Structural and Multidisciplinary
  Optimization 51~(3) (2015) 631--643.

\bibitem{feppon:hal-01972915}
F.~Feppon, G.~Allaire, C.~Dapogny,
  \href{https://hal.archives-ouvertes.fr/hal-01972915}{{Null space gradient
  flows for constrained optimization with applications to shape optimization}},
  working paper or preprint (Jan. 2019).
\newline\urlprefix\url{https://hal.archives-ouvertes.fr/hal-01972915}

\bibitem{allaire2011topology}
G.~Allaire, C.~Dapogny, P.~Frey, Topology and geometry optimization of elastic
  structures by exact deformation of simplicial mesh, Comptes Rendus
  Mathematique 349~(17-18) (2011) 999--1003.

\bibitem{yamasaki2011level}
S.~Yamasaki, T.~Nomura, A.~Kawamoto, K.~Sato, S.~Nishiwaki, A level set-based
  topology optimization method targeting metallic waveguide design problems,
  International Journal for Numerical Methods in Engineering 87~(9) (2011)
  844--868.

\bibitem{burman2018shape}
E.~Burman, D.~Elfverson, P.~Hansbo, M.~G. Larson, K.~Larsson, Shape
  optimization using the cut finite element method, Computer Methods in Applied
  Mechanics and Engineering 328 (2018) 242--261.

\bibitem{villanueva2017cutfem}
C.~H. Villanueva, K.~Maute, Cutfem topology optimization of 3d laminar
  incompressible flow problems, Computer Methods in Applied Mechanics and
  Engineering 320 (2017) 444--473.

\bibitem{belytschko2003topology}
T.~Belytschko, S.~Xiao, C.~Parimi, Topology optimization with implicit
  functions and regularization, International Journal for Numerical Methods in
  Engineering 57~(8) (2003) 1177--1196.

\bibitem{kreissl2012levelset}
S.~Kreissl, K.~Maute, Levelset based fluid topology optimization using the
  extended finite element method, Structural and Multidisciplinary Optimization
  46~(3) (2012) 311--326.

\bibitem{allaire2004structural}
G.~Allaire, F.~Jouve, A.-M. Toader, Structural optimization using sensitivity
  analysis and a level-set method, Journal of computational physics 194~(1)
  (2004) 363--393.

\bibitem{abhyankar2018petscts}
S.~Abhyankar, J.~Brown, E.~M. Constantinescu, D.~Ghosh, B.~F. Smith, H.~Zhang,
  Petsc/ts: A modern scalable ode/dae solver library (2018).
\newblock \href {http://arxiv.org/abs/1806.01437} {\path{arXiv:1806.01437}}.

\bibitem{adamsHighorderEllipticPDE2019}
T.~Adams, S.~Giani, W.~M. Coombs,
  \href{https://linkinghub.elsevier.com/retrieve/pii/S0021999118307915}{A
  high-order elliptic {{PDE}} based level set reinitialisation method using a
  discontinuous {{Galerkin}} discretisation}, Journal of Computational Physics
  379 (2019) 373--391.
\newblock \href {https://doi.org/10.1016/j.jcp.2018.12.003}
  {\path{doi:10.1016/j.jcp.2018.12.003}}.
\newline\urlprefix\url{https://linkinghub.elsevier.com/retrieve/pii/S0021999118307915}

\bibitem{geuzaine2009gmsh}
C.~Geuzaine, J.-F. Remacle, Gmsh: A 3-d finite element mesh generator with
  built-in pre-and post-processing facilities, International journal for
  numerical methods in engineering 79~(11) (2009) 1309--1331.

\bibitem{Rathgeber2016}
F.~Rathgeber, D.~A. Ham, L.~Mitchell, M.~Lange, F.~Luporini, A.~T.~T. McRae,
  G.-T. Bercea, G.~R. Markall, P.~H.~J. Kelly,
  \href{http://arxiv.org/abs/1501.01809}{Firedrake: automating the finite
  element method by composing abstractions}, ACM Trans. Math. Softw. 43~(3)
  (2016) 24:1--24:27.
\newblock \href {http://arxiv.org/abs/1501.01809} {\path{arXiv:1501.01809}},
  \href {https://doi.org/10.1145/2998441} {\path{doi:10.1145/2998441}}.
\newline\urlprefix\url{http://arxiv.org/abs/1501.01809}

\bibitem{Luporini2016}
F.~Luporini, D.~A. Ham, P.~H.~J. Kelly,
  \href{http://arxiv.org/abs/1604.05872}{An algorithm for the optimization of
  finite element integration loops}, ACM Transactions on Mathematical Software
  44 (2017) 3:1--3:26.
\newblock \href {http://arxiv.org/abs/1604.05872} {\path{arXiv:1604.05872}},
  \href {https://doi.org/10.1145/3054944} {\path{doi:10.1145/3054944}}.
\newline\urlprefix\url{http://arxiv.org/abs/1604.05872}

\bibitem{Homolya2017}
M.~Homolya, L.~Mitchell, F.~Luporini, D.~A. Ham,
  \href{http://arxiv.org/abs/1705.003667}{{TSFC: a structure-preserving form
  compiler}} (2017).
\newblock \href {http://arxiv.org/abs/1705.03667} {\path{arXiv:1705.03667}}.
\newline\urlprefix\url{http://arxiv.org/abs/1705.003667}

\bibitem{petsc-efficient}
S.~Balay, W.~D. Gropp, L.~C. McInnes, B.~F. Smith, Efficient management of
  parallelism in object oriented numerical software libraries, in: E.~Arge,
  A.~M. Bruaset, H.~P. Langtangen (Eds.), Modern Software Tools in Scientific
  Computing, Birkh{\"{a}}user Press, 1997, pp. 163--202.

\bibitem{petsc-user-ref}
S.~Balay, S.~Abhyankar, M.~F. Adams, J.~Brown, P.~Brune, K.~Buschelman,
  L.~Dalcin, V.~Eijkhout, W.~D. Gropp, D.~Karpeyev, D.~Kaushik, M.~G. Knepley,
  D.~A. May, L.~C. McInnes, R.~T. Mills, T.~Munson, K.~Rupp, P.~Sanan, B.~F.
  Smith, S.~Zampini, H.~Zhang, H.~Zhang, {PETS}c users manual, Tech. Rep.
  ANL-95/11 - Revision 3.11, Argonne National Laboratory (2019).

\bibitem{Dalcin2011}
L.~D. Dalcin, R.~R. Paz, P.~A. Kler, A.~Cosimo, Parallel distributed computing
  using {P}ython, Advances in Water Resources 34~(9) (2011) 1124--1139, new
  Computational Methods and Software Tools.
\newblock \href
  {https://doi.org/http://dx.doi.org/10.1016/j.advwatres.2011.04.013}
  {\path{doi:http://dx.doi.org/10.1016/j.advwatres.2011.04.013}}.

\bibitem{patankarCalculationProcedureHeat1972}
S.~Patankar, D.~Spalding,
  \href{https://linkinghub.elsevier.com/retrieve/pii/0017931072900543}{A
  calculation procedure for heat, mass and momentum transfer in
  three-dimensional parabolic flows}, International Journal of Heat and Mass
  Transfer 15~(10) (1972) 1787--1806.
\newblock \href {https://doi.org/10.1016/0017-9310(72)90054-3}
  {\path{doi:10.1016/0017-9310(72)90054-3}}.
\newline\urlprefix\url{https://linkinghub.elsevier.com/retrieve/pii/0017931072900543}

\bibitem{saadfgmres}
Y.~Saad, \href{https://doi.org/10.1137/0914028}{A flexible inner-outer
  preconditioned gmres algorithm}, SIAM J. Sci. Comput. 14~(2) (1993)
  461–469.
\newblock \href {https://doi.org/10.1137/0914028} {\path{doi:10.1137/0914028}}.
\newline\urlprefix\url{https://doi.org/10.1137/0914028}

\bibitem{burman2015cutfem}
E.~Burman, S.~Claus, P.~Hansbo, M.~G. Larson, A.~Massing, Cutfem: discretizing
  geometry and partial differential equations, International Journal for
  Numerical Methods in Engineering 104~(7) (2015) 472--501.

\bibitem{belytschko2003structured}
T.~Belytschko, C.~Parimi, N.~Mo{\"e}s, N.~Sukumar, S.~Usui, Structured extended
  finite element methods for solids defined by implicit surfaces, International
  journal for numerical methods in engineering 56~(4) (2003) 609--635.

\bibitem{zenodo/Firedrake-20211025.1}
\href{https://doi.org/10.5281/zenodo.5596324}{{Software used in
  `Three-dimensional topology optimization of heat exchangers with the
  level-set method'}} (oct 2021).
\newblock \href {https://doi.org/10.5281/zenodo.5596324}
  {\path{doi:10.5281/zenodo.5596324}}.
\newline\urlprefix\url{https://doi.org/10.5281/zenodo.5596324}

\bibitem{miguel_salazar_de_troya_2021_5597909}
M.~Salazar~de Troya,
  \href{https://doi.org/10.5281/zenodo.5597909}{{LLNL}/lestofire: First
  release} (Oct. 2021).
\newblock \href {https://doi.org/10.5281/zenodo.5597909}
  {\path{doi:10.5281/zenodo.5597909}}.
\newline\urlprefix\url{https://doi.org/10.5281/zenodo.5597909}

\bibitem{paroliniComputationalFluidDynamics2005}
N.~Parolini, \href{http://infoscience.epfl.ch/record/33607}{Computational fluid
  dynamics for naval engineering problems}, Journal of Computational Physics
  (2005).
\newblock \href {https://doi.org/10.5075/EPFL-THESIS-3138}
  {\path{doi:10.5075/EPFL-THESIS-3138}}.
\newline\urlprefix\url{http://infoscience.epfl.ch/record/33607}

\end{thebibliography}

\appendix
\section{Signed distance equation solver}
\label{sec:sds}
We cannot apply the Dirichlet boundary conditions in Equation \eqref{eq:signed_distance} because the isocontour $\phi = 0$ does not typically match the mesh elements.
To overcome this issue, we follow \cite{paroliniComputationalFluidDynamics2005} and first calculate the local $L^2$ projection $\varphi_{\text {int }} \in H^1(D_{\text{int}})$ such that
\begin{align}
	\int_{D_{\text {int }}} \varphi_{\text {int }} v ~dV=\int_{D_{\text {int }}} \frac{\phi}{|\nabla \phi|} v ~dV\,,
\end{align}
for all $v \in H^1(D_{\text{int}})$.
We then apply $\varphi_{\text{int}}$ as a Dirichlet condition in the subdomain $D_{\text {int }}=\left\{\bigcup K: K \in \mathcal{T}_{h}, K \cap \Gamma \neq \varnothing\right\}$, i.e. all the elements crossed by the isocontour $\phi = 0$.

To solve the nonlinear elliptic Equation \eqref{eq:signed_distance}, we apply the Picard iteration method:
Find the iterate $\varphi^{(m)} \in W_{\phi} = \{v_{h} \in H^{1}(D):\left.v_{h}\right|_{\mathcal{K}} \in P^{1}(\mathcal{K})~ \forall \mathcal{K} \in \mathcal{T}_{h} ~|~ v_{h} = \varphi_{\text{int}} ~\text{in}~  D_{\text {int }}\}$
\begin{equation}
	\int_{D}\nabla \varphi^{(m)} \cdot \nabla v~dV=  \int_{D}\left( 1 - \iota(| \nabla \varphi^{(m-1)}|) \right) \nabla \varphi^{(m-1)} \cdot \nabla v ~dV
\end{equation}
$\forall ~v\in W_{\phi}$, update $\varphi^{(m-1)} = \varphi^{(m)}$ and $m \leftarrow m+1$, repeat until convergence and update the level-set, i.e. $\phi = \varphi^{(m)}$.
In practice, convergence is attained in less than 10 iterations.

\end{document}